\documentclass[10pt]{article}

\usepackage{amsmath}
\usepackage{amsthm}
\usepackage{amssymb}
\usepackage{amscd}
\usepackage{amsfonts}  
\usepackage{makeidx}
\usepackage{hyperref}
\usepackage{makeidx}     
\usepackage{graphicx}    
\usepackage{multicol}    

\DeclareMathAlphabet\EuR{U}{eur}{m}{n}
\SetMathAlphabet\EuR{bold}{U}{eur}{b}{n}

\begin{document}


\newtheorem{theorem}{Theorem}[section]
\newtheorem{lemma}[theorem]{Lemma}
\newtheorem{proposition}[theorem]{Proposition}
\newtheorem{definition}[theorem]{Definition}
\newtheorem{example}[theorem]{Example}
\newtheorem{remark}[theorem]{Remark}
\newtheorem{corollary}[theorem]{Corollary}
\newtheorem{conjecture}[theorem]{Conjecture}
\newtheorem{problem}[theorem]{Problem}
{\catcode`@=11\global\let\c@equation=\c@theorem}
\renewcommand{\theequation}{\thetheorem}

\renewcommand{\theenumi}{\alph{enumi}}
\renewcommand{\labelenumi}{(\theenumi)}

\makeatletter
\renewcommand{\@seccntformat}[1]{\csname the#1\endcsname.\hspace{1em}}
\makeatother
\newcommand{\tit}[1]{\begin{bf} \begin{center} \begin{Large}
\section{#1}
\label{sec: #1}
\end{Large}\end{center}\end{bf}
\nopagebreak}


\newcommand{\squarematrix}[4]{\left( \begin{array}{cc} #1 & #2 \\ #3 &
#4
\end{array} \right)}
\newcommand{\smallmx}[4]{\mbox{\begin{scriptsize}$\squarematrix{#1}{#2}
        {#3}{#4}$\end{scriptsize}}}

\newcommand{\comsquare}[8]{
\begin{center}
$\begin{CD}
#1 @>#2>> #3\\
@V{#4}VV @VV{#5}V\\
#6 @>>#7> #8
\end{CD}$
\end{center}}


\let\sect=\S
\newcommand{\curs}{\EuR}
\newcommand{\CHAINCOMPLEXES}{\curs{CHCOM}}
\newcommand{\GROUPOIDS}{\curs{GROUPOIDS}}
\newcommand{\PAIRS}{\curs{PAIRS}}
\newcommand{\FGINJ}{\curs{FGINJ}}
\newcommand{\Or}{\curs{Or}}
\newcommand{\SPACES}{\curs{SPACES}}
\newcommand{\SPECTRA}{\curs{SPECTRA}}
\newcommand{\Sub}{\curs{Sub}}


\newcommand{\bbC}{{\mathbb C}}
\newcommand{\bbH}{{\mathbb H}}
\newcommand{\bbI}{{\mathbb I}}
\newcommand{\bbK}{{\mathbb K}}
\newcommand{\bbKO}{\mathbb{KO}}
\newcommand{\bbN}{{\mathbb N}}
\newcommand{\bbP}{{\mathbb P}}
\newcommand{\bbQ}{{\mathbb Q}}
\newcommand{\bbR}{{\mathbb R}}
\newcommand{\bbZ}{{\mathbb Z}}

\newcommand{\calbh}{{\mathcal B}{\mathcal H}}
\newcommand{\calc}{{\mathcal C}}
\newcommand{\cald}{{\mathcal D}}
\newcommand{\cale}{{\mathcal E}}
\newcommand{\calf}{{\mathcal F}}
\newcommand{\calg}{{\mathcal G}}
\newcommand{\calh}{{\mathcal H}}
\newcommand{\calk}{{\mathcal K}}

\newcommand{\bfA}{\ensuremath{\mathbf{A}}}
\newcommand{\bfa}{\ensuremath{\mathbf{a}}}
\newcommand{\bfb}{\ensuremath{\mathbf{b}}}
\newcommand{\bfCTR}{\ensuremath{\mathbf{CTR}}}
\newcommand{\bfDtr}{\ensuremath{\mathbf{Dtr}}}
\newcommand{\bfE}{\ensuremath{\mathbf{E}}}
\newcommand{\bff}{\ensuremath{\mathbf{f}}}
\newcommand{\bfF}{\ensuremath{\mathbf{F}}}
\newcommand{\bfg}{\ensuremath{\mathbf{g}}}
\newcommand{\bfHH}{\ensuremath{\mathbf{HH}}}
\newcommand{\bfHPB}{\ensuremath{\mathbf{HPB}}}
\newcommand{\bfI}{\ensuremath{\mathbf{I}}}
\newcommand{\bfi}{\ensuremath{\mathbf{i}}}
\newcommand{\bfK}{\ensuremath{\mathbf{K}}}
\newcommand{\bfL}{\ensuremath{\mathbf{L}}}
\newcommand{\bfr}{\ensuremath{\mathbf{r}}}
\newcommand{\bfs}{\ensuremath{\mathbf{s}}}
\newcommand{\bfS}{\ensuremath{\mathbf{S}}}
\newcommand{\bft}{\ensuremath{\mathbf{t}}}
\newcommand{\bfT}{\ensuremath{\mathbf{T}}}
\newcommand{\bfTC}{\ensuremath{\mathbf{TC}}}
\newcommand{\bfU}{\ensuremath{\mathbf{U}}}
\newcommand{\bfu}{\ensuremath{\mathbf{u}}}
\newcommand{\bfv}{\ensuremath{\mathbf{v}}}
\newcommand{\bfw}{\ensuremath{\mathbf{w}}}


\newcommand{\aut}{\operatorname{aut}}
\newcommand{\Bor}{\operatorname{Bor}}
\newcommand{\ch}{\operatorname{ch}}
\newcommand{\class}{\operatorname{class}}
\newcommand{\cok}{\operatorname{coker}}
\newcommand{\cone}{\operatorname{cone}}
\newcommand{\colim}{\operatorname{colim}}
\newcommand{\con}{\operatorname{con}}
\newcommand{\conhom}{\operatorname{conhom}}
\newcommand{\cyclic}{\operatorname{cyclic}}
\newcommand{\ev}{\operatorname{ev}}
\newcommand{\ext}{\operatorname{ext}}
\newcommand{\Gen}{\operatorname{Gen}}
\newcommand{\hur}{\operatorname{hur}}
\newcommand{\im}{\operatorname{im}}
\newcommand{\inj}{\operatorname{inj}}
\newcommand{\id}{\operatorname{id}}
\newcommand{\infl}{\operatorname{Infl}}
\newcommand{\ind}{\operatorname{ind}}
\newcommand{\Inn}{\operatorname{Inn}}
\newcommand{\Irr}{\operatorname{Irr}}
\newcommand{\Is}{\operatorname{Is}}
\newcommand{\map}{\operatorname{map}}
\newcommand{\MOD}{\operatorname{MOD}}
\newcommand{\mor}{\operatorname{mor}}
\newcommand{\Ob}{\operatorname{Ob}}
\newcommand{\op}{\operatorname{op}}
\newcommand{\pr}{\operatorname{pr}}
\newcommand{\point}{\operatorname{pt.}}
\newcommand{\Rat}{\operatorname{Rat}}
\newcommand{\res}{\operatorname{res}}
\newcommand{\topo}{\operatorname{top}}
\newcommand{\tors}{\operatorname{tors}}

\newcommand{\pt}{\{\point \}}


\newcommand{\comment}[1]                      
{
{{\bf Comment: }{\ttfamily #1}}
}

\newcommand{\tabtit}[1]
{
\ref{sec: #1}. & #1
}


\title{Equivariant Cohomological Chern Characters}
\author{
Wolfgang L\"uck\thanks{\noindent email:
lueck@math.uni-muenster.de\protect\\
www: ~http://www.math.uni-muenster.de/u/lueck/\protect\\
FAX: 49 251 8338370\protect}\\
Fachbereich Mathematik\\ Universit\"at M\"unster\\
Einsteinstr.~62\\ 48149 M\"unster\\Germany}
\maketitle

\typeout{-----------------------  Abstract  ------------------------}
\begin{abstract}
We construct for an equivariant cohomology theory
for proper equivariant $CW$-complexes an equivariant Chern character, provided
that certain conditions about the coefficients are satisfied. These conditions are fulfilled
if the coefficients of the equivariant cohomology theory possess a
Mackey structure. Such a structure is present in many interesting examples.

\smallskip

\noindent
Key words: equivariant cohomology theory, equivariant Chern character. \\
Mathematics Subject Classification 2000: 55N91.
\end{abstract}

\typeout{-----------------------  Introduction ------------------------}

\setcounter{section}{-1}
\tit{Introduction}

The purpose of this paper is to construct an equivariant Chern character
for a proper equivariant cohomology theory $\calh^*_?$ with values in $R$-modules 
for a commutative associative ring $R$ with unit which satisfies 
$\bbQ \subseteq R$. It is a natural transformation of
equivariant cohomology theories
$$\ch^*_? \colon \calh^*_? ~ \to \calbh^*_?$$
for a given equivariant cohomology theory $\calh^*_?$.
Here $\calbh^*_?$ is the associated equivariant cohomology theory
which is defined by the Bredon cohomology with coefficients coming from 
the coefficients of $\calh^*_?$. The notion of an equivariant cohomology theory and
examples for it are presented in Section~\ref{sec: Equivariant Cohomology Theories} and the 
associated Bredon cohomology is explained in 
Section~\ref{sec: The Associated Bredon Cohomology Theory}. 
The point is that $\calbh^*_?$ is much simpler and easier
to compute than $\calh^*_?$. If $\calh^*_?$ satisfies the disjoint union axiom,
then 
$$\ch^n_G(X,A) \colon \calh^n_G(X,A) \xrightarrow{\cong} \calbh^n_G(X,A)$$ 
is bijective for every discrete group $G$, proper $G$-$CW$-pair
$(X,A)$ and $n \in \bbZ$. 

The Chern character $\ch^*_?$ is only defined
if the coefficients of $\calh_?^*$ satisfy a certain injectivity condition
(see Theorem~\ref{the: equivariant Chern character and injective coeff}).
This condition is fulfilled if the coefficients of $\calh^*_?$
come with a Mackey structure 
(see Theorem~\ref{the: Chern for calh^*_G with a Mackey structure}) 
what is the case in many interesting examples.

The equivariant cohomological Chern character is a generalization to the equivariant setting
of the classical non-equivariant Chern character for a (non-equivariant) cohomology theory 
$\calh^*$ (see Example~\ref{exa: cohomological version of Dold's construction})
$$ \ch^n(X,A) \colon \calh^n(X,A)  \xrightarrow{\cong} \prod_{p+q = n} H^p(X,A,\calh^q(*)).
$$
The equivariant cohomological Chern character has already been constructed in the special case, where
$\calh^*_?$ is equivariant topological $K$-theory $K^*_?$, in
\cite{Lueck-Oliver(2001b)}. Its homological
version has already been treated in \cite{Lueck(2002b)}
and plays an important role in the computation of the source of the assembly
maps appearing in the Farrell-Jones Conjecture for
$K_n(RG)$ and $L^{\langle - \infty \rangle}_n(RG)$ and the Baum-Connes Conjecture
for $K_n(C^*_r(G))$ (see also \cite{Lueck(2002d)}).

The detailed formulation of the  main result of this paper is
presented in Theorem~\ref{the: Chern for calh^*_G with a Mackey structure}.

The equivariant Chern character will play a key role in the proof of the following
result which will be presented in \cite{Lueck(2004i)}.

\begin{theorem}[Rational computation of the topological $K$-theory of $BG$]
 \label{the: rational computation of K^*(BG)}
Let $G$ be a discrete group. Suppose that there is a 
finite $G$-$CW$-model for the classifying space $\underline{E}G$ for
proper $G$-actions.
Then there is a $\bbQ$-isomorphism,
natural in $G$ and compatible with the multiplicative structures
\begin{multline*}
\overline{\ch}^n_G \colon K^n(BG) \otimes_{\bbZ} \bbQ
\\ \xrightarrow{\cong} ~
\prod_{i \in \bbZ} H^{2i+n}(BG;\bbQ) \times
\prod_{p \text{ prime}} ~ \prod_{(g) \in \con_p(G)} 
H^{2i+n}(BC_G\langle g \rangle;\bbQ\widehat{_p}).
\end{multline*}
\end{theorem}

Here $\con_p(G)$ is the set of conjugacy classes (g) of elements $g
\in G$ of order $p^m$ for some integer $m \ge 1$ and $C_G\langle g
\rangle$ is the centralizer of
the cyclic subgroup generated by $g$ in $G$.

The assumption in Theorem~\ref{the: rational computation of K^*(BG)} that there is a 
finite $G$-$CW$-model for the classifying space $\underline{E}G$ for
proper $G$-actions is satisfied
for instance, if $G$ is word-hyperbolic in the sense of Gromov,
if $G$ is a cocompact subgroup of a Lie group with finitely many path
components, if $G$ is a finitely generated one-relator group,
if $G$ is an arithmetic group, a mapping class group of a compact surface or the group
of outer automorphisms of a finitely generated free group. For more information about
$\underline{E}G$ we refer for instance to
\cite{Baum-Connes-Higson(1994)} and \cite{Lueck(2004h)}. 

A group $G$ is always understood  to be discrete and a ring $R$ is always understood to be
associative with unit throughout this paper.

The paper is organized as follows:

\begin{tabular}{ll}
\tabtit{Equivariant Cohomology Theories}
\\
\tabtit{Modules over a Category}
\\
\tabtit{The Associated Bredon Cohomology Theory}
\\
\tabtit{The Construction of the Equivariant Cohomological Chern Character}
\\
\tabtit{Mackey Functors}
\\
\tabtit{Multiplicative Structures}
\\
 & References
\end{tabular}
\bigskip


\setcounter{section}{0}

\typeout{-----------------------  Section 1  ------------------------}

\tit{Equivariant Cohomology Theories}

In this section we describe the axioms
of a (proper) equivariant cohomology theory. They are dual to the ones
of a (proper) equivariant homology theory as described in
\cite[Section 1]{Lueck(2002b)}.

Fix a group $G$ and an commutative ring $R$. 
A $G$-$CW$-pair $(X,A)$ is a pair of $G$-$CW$-complexes.
It is \emph{proper} if all isotropy groups of $X$ are finite.
It is \emph{relative finite} if $X$ is obtained from $A$ by attaching
finitely many equivariant cells, or, equivalently, if
$G\backslash (X/A)$ is compact.
Basic information about $G$-$CW$-pairs can be found for instance in
\cite[Section 1 and 2]{Lueck(1989)}.
A \emph{$G$-cohomology theory $\calh^*_G$ with values
in $R$-modules} is a collection of
covariant functors $\calh_G^n$ from the category of
$G$-$CW$-pairs to the category of
$R$-modules indexed by $n \in \bbZ$ together with natural transformations
$\delta^n_G(X,A)\colon \calh^n_G(X,A) \to
\calh^{n+1}_G(A):= \calh^{n+1}_G(A,\emptyset)$ for $n \in \bbZ$
such that the following axioms are satisfied:
\begin{itemize}

\item $G$-homotopy invariance \\[1mm]
If $f_0$ and $f_1$ are $G$-homotopic maps $(X,A) \to (Y,B)$
of $G$-$CW$-pairs, then $\calh^n_G(f_0) = \calh^n_G(f_1)$ for $n \in \bbZ$;

\item Long exact sequence of a pair\\[1mm]
Given a pair $(X,A)$ of $G$-$CW$-complexes,
there is a long exact sequence
$$\ldots \xrightarrow{\delta^{n-1}_G}
\calh^{n}_G(X,A) \xrightarrow{\calh^{n}_G(j)} 
\calh^n_G(X) \xrightarrow{\calh^{n}_G(i)} \calh^n_G(A)
\xrightarrow{\delta^n_G} \ldots,$$
where $i\colon A \to X$ and $j\colon X \to (X,A)$ are the inclusions;

\item Excision \\[1mm]
Let $(X,A)$ be a $G$-$CW$-pair and let
$f\colon A \to B$ be a cellular $G$-map of
$G$-$CW$-complexes. Equip $(X\cup_f B,B)$ with the induced structure
of a $G$-$CW$-pair. Then the canonical map
$(F,f)\colon (X,A) \to (X\cup_f B,B)$ induces an isomorphism
$$\calh^n_G(F,f)\colon \calh^n_G(X,A) \xrightarrow{\cong}
\calh^n_G(X\cup_f B,B).$$

\end{itemize}

Sometimes also the following axiom is required.

\begin{itemize}

\item Disjoint union axiom\\[1mm]
Let $\{X_i\mid i \in I\}$ be a family of
$G$-$CW$-complexes. Denote by
$j_i\colon  X_i \to \coprod_{i \in I} X_i$ the canonical inclusion.
Then the map
$$\prod_{i \in I} \calh^{n}_G(j_i)\colon  \calh^n_G\left(\coprod_{i \in I} X_i\right)
\xrightarrow{\cong}  \prod _{i \in I} \calh^n_G(X_i)
$$
is bijective.

\end{itemize}

If $\calh^*_G$ is defined or considered only for proper $G$-$CW$-pairs
$(X,A)$, we call it a
\emph{proper $G$-cohomology theory $\calh^*_G$ with values
in $R$-modules}.

The role of the disjoint union axiom is explained by the following result.
Its proof for non-equivariant cohomology theories 
(see for instance \cite[7.66 and 7.67]{Switzer(1975)}) 
carries over directly  to $G$-cohomology theories.

\begin{lemma} \label{lem: role of the disjoint union axiom}
Let $\calh^*_G$ and $\calk^*_G$ be (proper) $G$-cohomology theories. Then

\begin{enumerate}

\item \label{lem: role of the disjoint union axiom: Milnor's sequence}
Suppose that $\calh^*_G$ satisfies the disjoint union axiom. Then there exists for every
(proper) $G$-$CW$-pair $(X,A)$ a natural short exact sequence
\begin{multline*}
0 \to {\lim}^1_{n \to \infty} \calh^{p-1}_G(X_n \cup A,A) \to \calh^p(X,A) \to 
\lim_{n \to \infty} \calh^p_G(X_n \cup A,A) \to 0;
\end{multline*}

\item \label{lem: role of the disjoint union axiom: transformations of equiv. coho. th.}
Let $T^* \colon \calh^*_G \to \calk^*_G$ be a transformation of (proper) $G$-cohomology theories,
i.e. a collection of  natural transformations $T^n \colon \calh^n_G \to \calk^n_G$
of contravariant functors from the category of (proper)
$G$-$CW$-pairs to the category of $R$-modules indexed by $n \in \bbZ$ which is compatible
with the boundary operator associated to (proper) $G$-$CW$-pairs. Suppose that
$T^n(G/H)$ is bijective for every (proper) homogeneous space $G/H$ and
$n \in \bbZ$.

Then  $T^n(X,A) \colon \calh^*_G(X,A) \to \calk^*_G(X,A)$ is bijective for all 
$n \in \bbZ$ provided that $(X,A)$ is relative finite or 
that both $\calh^*$ and $\calk^*$ satisfy the disjoint union axiom.

\end{enumerate}
\end{lemma}

\begin{remark}[The disjoint union axiom is not compatible with $- \otimes_{\bbZ} \bbQ$]
\label{rem: The disjoint union axiom is not compatible with - otimes_{bbZ} bbQ} \em
Let $\calh^*_G$ be a $G$-cohomology theory with values in $\bbZ$-modules.
Then $\calh^*_G \otimes_{\bbZ} \bbQ $ is a $G$-cohomology theory with values in
$\bbQ$-modules since $\bbQ$ is flat as $\bbZ$-module. However,
even if $\calh^*$ satisfies the disjoint union axiom, 
$\calh^*_G \otimes_{\bbZ} \bbQ$ does not satisfy the disjoint union axiom
since $- \otimes_{\bbZ} \bbQ$ is not compatible with products over arbitrary index sets.
\em
\end{remark}

\begin{example}[Rationalizing topological $K$-theory] 
\label{exa: Rationalizing topological  K -theory} \em
Consider for instance the (non-equivariant) cohomology theory with values in $\bbZ$-modules
satisfying the disjoint union axiom given by topological $K$-theory $K^*$. 
Let $K^*(-;\bbQ)$ be the cohomology theory associated to the rationalization
of the $K$-theory spectrum. This is a (non-equivariant) cohomology theory with values 
in $\bbQ$-modules satisfying the disjoint union axiom. There is a natural transformation
$$T^*(X) \colon K^*(X) \otimes_{\bbZ} \bbQ ~ \to K^*(X;\bbQ)$$
of (non-equivariant) cohomology theory with values in $\bbQ$-modules. 
The $\bbQ$-map $T^n(\pt)$ is bijective for all $n \in \bbZ$. Hence
$T^n(X)$ is bijective for all finite $CW$-complexes by 
Lemma \ref{lem: role of the disjoint union axiom}
\eqref{lem: role of the disjoint union axiom: transformations of equiv. coho. th.}.
Notice that $K^*(X) \otimes_{\bbZ} \bbQ$ does not satisfy the disjoint union axiom
in contrast to $K^*(X;\bbQ)$. Hence we cannot expect $T^n(X)$ to be bijective for all
$CW$-complexes. Consider the case $X = BG$ for a  finite group $G$.
Since $H^p(BG;\bbQ) \cong H^p(\pt;\bbQ)$ for all $p \in \bbZ$, one obtains
$K^p(BG;\bbQ) \cong K^p(\pt;\bbQ)$ for all $p \in \bbZ$. By the Atiyah-Segal Completion
Theorem $K^p(BG) \otimes_{\bbZ} \bbQ \cong K^p(\pt) \otimes_{\bbZ} \bbQ$ is only true
if and only if the finite group $G$ is trivial. 
\em
\end{example}

Let $\alpha\colon H \to G$ be a group homomorphism.
Given an $H$-space $X$, define the \emph{induction of $X$ with $\alpha$}
to be the $G$-space $\ind_{\alpha} X$ which  is the quotient of
$G \times X$ by the right $H$-action
$(g,x) \cdot h := (g\alpha(h),h^{-1} x)$
for $h \in H$ and $(g,x) \in G \times X$.
If $\alpha\colon H \to G$ is an inclusion, we also write $\ind_H^G$ instead of
$\ind_{\alpha}$.

A \emph{(proper) equivariant cohomology theory
$\calh^*_?$ with values in $R$-modules}
consists of a collection of (proper)
$G$-cohomology theory $\calh^*_G$ with values in $R$-modules for each group $G$
together with the following so called \emph{induction structure}:
given a group homomorphism $\alpha\colon  H \to G$ and  a (proper) $H$-$CW$-pair
$(X,A)$ such that $\ker(\alpha)$ acts freely on $X$,
there are for each $n \in \bbZ$
natural isomorphisms
\begin{eqnarray}
\ind_{\alpha}\colon \calh^n_G(\ind_{\alpha}(X,A)) 
&\xrightarrow{\cong} &
\calh^n_H(X,A)\label{induction structure}
\end{eqnarray}
satisfying 

\begin{enumerate}

\item Compatibility with the boundary homomorphisms\\[1mm]
$\delta^n_H \circ \ind_{\alpha} = \ind_{\alpha} \circ \delta^n_G$;

\item Functoriality\\[1mm]
Let $\beta\colon  G \to K$ be another group
homomorphism such that $\ker(\beta \circ \alpha)$ acts freely on $X$.
Then we have for $n \in \bbZ$
$$\ind_{\beta \circ \alpha} ~ = ~
\ind_{\alpha} \circ \ind_{\beta} \circ \calh^n_K(f_1) \colon
\calh^n_H(\ind_{\beta\circ\alpha}(X,A)) \to \calh^n_K(X,A),$$
where $f_1\colon \ind_{\beta}\ind_{\alpha}(X,A)
\xrightarrow{\cong} \ind_{\beta\circ \alpha}(X,A),
\hspace{1mm} (k,g,x) \mapsto (k\beta(g),x)$
is the natural $K$-homeo\-mor\-phism;

\item Compatibility with conjugation\\[1mm]
For $n \in \bbZ$, $g \in G$ and a (proper) $G$-$CW$-pair $(X,A)$
the homomorphism $\ind_{c(g)\colon G \to G}\colon 
\calh^n_G(\ind_{c(g)\colon G \to G}(X,A)) \to \calh^n_G(X,A)$ agrees with
$\calh^n_G(f_2)$, where $f_2$ is the $G$-homeomorphism
$f_2\colon (X,A) \to \ind_{c(g)\colon  G \to G} (X,A), \hspace{1mm} 
x \mapsto(1,g^{-1}x)$ and $c(g)(g') = gg'g^{-1}$.

\end{enumerate}

This induction structure links the
various $G$-cohomology theories for different groups $G$.
It will play a key role in the construction
of the equivariant  Chern character even if we want to carry it out only
for a fixed group $G$.
In all of the relevant examples the induction homomorphism $\ind_{\alpha}$
of \eqref{induction structure}
exists for every group homomorphism $\alpha\colon H \to G$, the condition
that $\ker(\alpha)$ acts freely on $X$ is only needed to ensure
that $\ind_{\alpha}$ is bijective. If $\alpha$ is an inclusion, we sometimes write
$\ind_H^G$ instead of $\ind_{\alpha}$. 

We say that $\calh^*_?$ satisfies the \emph{disjoint union axiom} if for every group $G$ 
the $G$-cohomology theory $\calh^*_G$ satisfies the disjoint union axiom.

We will later need the following lemma 
whose elementary proof is analogous to the 
one in \cite[Lemma 1.2]{Lueck(2002b)}.
\begin{lemma}
\label{lem: calh_G(G/H) and calh_H(*)}
Consider finite subgroups $H,K \subseteq G$ and an element
$g \in G$ with $gHg^{-1} \subseteq K$.
Let $R_{g^{-1}}\colon G/H \to G/K$ be the $G$-map
sending $g'H$ to $g'g^{-1}K$ and
$c(g)\colon H \to K$ be the homomorphism sending $h$ to $ghg^{-1}$.
Let $\pr\colon (\ind_{c(g)\colon H \to K}\pt) \to \pt$ be the projection.
Then the following diagram commutes
$$\begin{CD}
\calh_G^n(G/K)   @> \calh^n_G(R_{g^{-1}}) >> \calh_G^n(G/H)
\\
@ V \ind_K^G V \cong V @V \ind_H^G V\cong V 
\\
\calh^n_K(\pt) @>  \ind_{c(g)} \circ \calh_K^n(\pr) >>  \calh^n_H(\pt) 
\end{CD}
$$
\end{lemma}

\begin{example}[Borel cohomology]
\label{exa: cohomology of the quotient or Borel construction}
\em
Let $\calk^*$ be a cohomology theory for (non-equivariant)
$CW$-pairs with values in $R$-modules.
Examples are singular cohomology and topological $K$-theory.
Then we obtain two equivariant cohomology
theories with values in $R$-modules by the following  constructions
\begin{eqnarray*}
\calh^n_G(X,A) & = & \calk^n(G\backslash X,G\backslash A);
\\
\calh^n_G(X,A) & = & \calk^n(EG \times_G (X,A)).
\end{eqnarray*}
The second one is called the
\emph{equivariant Borel cohomology associated to $\calk$}.
In both cases $\calh^*_G$ inherits the structure of a
$G$-cohomology theory from the cohomology structure on $\calk^*$.

The induction homomorphism associated to a group homomorphism
$\alpha\colon H \to G$ is defined as follows.
Let $a\colon  H\backslash X \xrightarrow{\cong} G\backslash (G \times_{\alpha} X)$
be the homeomorphism sending $Hx$ to $G(1,x)$.
Define $b\colon EH \times_H X \to EG \times_G G \times_{\alpha} X$
by sending $(e,x)$ to $(E\alpha(e),1,x)$ for $e \in EH$,
$x \in X$ and $E\alpha\colon EH \to EG$ the
$\alpha$-equivariant map induced by $\alpha$.
The desired induction map $\ind_{\alpha}$ is given
by $\calk^*(a)$ and $\calk^*(b)$. If the kernel $\ker(\alpha)$ acts freely on
$X$, the map $b$ is a homotopy equivalence and hence in both cases
$\ind_{\alpha}$ is bijective. 

If $\calk^*$ satisfies the disjoint union axiom, the same is true for the
two equivariant cohomology theories constructed above.
\em
\end{example}

\begin{example}[Equivariant $K$-theory] \label{exa: equivariant K-theory}  \em
In \cite{Lueck-Oliver(2001b)}  $G$-equivariant topological (complex) 
$K$-theory $K^*_G(X,A)$ is constructed for any proper
$G$-$CW$-pair $(X,A)$ and shown that $K^*_G$ defines a proper 
$G$-cohomology theory satisfying the disjoint union axiom. Given a group homomorphism
$\alpha \colon H \to G$, it induces an injective group homomorphism
$\overline{\alpha} \colon H/\ker(\alpha) \to G$. Let
$$\infl^H_{H/\ker(\alpha)} \colon
K_{H/\ker(\alpha)}^*(\ker(\alpha)\backslash X) \to K_H^*(X)$$ 
be the inflation homomorphism of \cite[Proposition 3.3]{Lueck-Oliver(2001b)}
and
$$\ind_{\overline{\alpha}} \colon
K_{H/\ker(\alpha)}^*(\ker(\alpha)\backslash X) \xrightarrow{\cong}
K_G^*(\ind_{\overline{\alpha}}(\ker(\alpha)\backslash X))$$
be the induction isomorphism of  
\cite[Proposition 3.2 (b)]{Lueck-Oliver(2001b)}. Define the induction homomorphism
$$\ind_{\alpha} \colon K^*_G(\ind_{\alpha} X) \to K^*_H(X)$$
by $\infl^H_{H/\ker(\alpha)} \circ
(\ind_{\overline{\alpha}})^{-1}$, where we identify
$\ind_{\alpha} X = \ind_{\overline{\alpha}}(\ker(\alpha)\backslash X)$.
On the level of complex finite-dimensional vector bundles the
induction homomorphism
$\ind_{\alpha}$ corresponds to considering for a $G$-vector bundle
$E$ over $G \times_{\alpha} X$ the $H$-vector bundle obtained from
$E$ by the pullback construction associated to the $\alpha$-equivariant map
$X \to G \times_{\alpha} X,~ x \mapsto (1,x)$. 

Thus we obtain a proper equivariant cohomology theory $K^*_?$ with values
in $\bbZ$-modules which satisfies the disjoint union axiom.
There is also a real version $KO^*_?$.
\em
\end{example}

\begin{example}[Equivariant cohomology theories and spectra] 
\label{exa: Equivariant cohomology theories and spectra} \em
Denote by $\GROUPOIDS$ the category of small groupoids. Let
$\Omega\text{-}\SPECTRA$ be the category of $\Omega$-spectra, where a 
morphism $\bff \colon \bfE \to \bfF$ is a sequence of maps
$f_n \colon E_n \to F_n$ compatible with the structure maps and 
we work in the category of compactly generated spaces
(see for instance \cite[Section 1]{Davis-Lueck(1998)}).
A contravariant $\GROUPOIDS$-$\Omega$-spectrum is a contravariant
functor $\bfE \colon \GROUPOIDS \to \Omega\text{-}\SPECTRA$. 

Next we explain how we can associate to it an equivariant
cohomology theory $H^*_?(-;\bfE)$ satisfying the disjoint union axiom, provided
that $\bfE$ respects equivalences, i.e.\ it sends an
equivalence of groupoids to a weak equivalence of spectra. 
This construction is dual to the
construction of an equivariant homology theory associated to a covariant
$\GROUPOIDS$-spectrum as explained in \cite[Section 6.2]{Lueck-Reich(2004g)},
\cite[Theorem~2.10 on page 21]{SauerJ(2002)}.

Fix a group $G$. We have to specify a $G$-cohomology theory $\calh^*_G(-;\bfE)$.
Let $\Or(G)$ be the \emph{orbit category}
whose set of objects consists of homogeneous $G$-spaces $G/H$ and whose morphisms
are $G$-maps. For a $G$-set $S$ we denote by $\calg^G(S)$ its associated 
\emph{transport  groupoid}. Its objects are the elements of $S$.
The set of morphisms from $s_0$ to $s_1$ consists
of those elements $g \in G$ which satisfy $g s_0 = s_1$. Composition in $\calg^G(S)$
comes from the multiplication in $G$. Thus we obtain for a group $G$ a covariant functor
\begin{eqnarray}
\calg^G\colon \Or(G) \to \GROUPOIDS, \hspace{7mm} G/H \mapsto \calg^G(G/H),
\label{functor calg^G}
\end{eqnarray}
and a contravariant $\Or(G)$-$\Omega$-spectrum $\bfE \circ \calg^G$. 
Given a $G$-$CW$-pair $(X,A)$, we obtain a contravariant pair of $\Or(G)$-$CW$-complexes
$(X^?,A^?)$
by sending $G/H$ to $(\map_G(G/H,X),\map_G(G/H,A)) = (X^H,A^H)$. 
The contravariant $\Or(G)$-spectrum $\bfE \circ \calg^G$ defines a cohomology theory
on the category of contravariant $\Or(G)$-$CW$-complexes as explained in 
\cite[Section 4]{Davis-Lueck(1998)}. It value at $(X^?,A^?)$ is defined to be $H^*_G(X,A;\bfE)$.
Explicitely, $H^n_G(X,A;\bfE)$ is the $(-n)$-th homotopy group of the 
spectrum $\map_{\Or(G)}\left(X_+^? \cup_{A_+^?} \cone(A_+^?),\bfE \circ \calg^G\right)$.
We need $\Omega$-spectra in order to ensure that the disjoint union axiom holds.

We briefly explain for a group homomorphism $\alpha \colon H \to G$
the definition of the induction homomorphism
$\ind_{\alpha} \colon \calh^n_G(\ind_{\alpha}X;\bfE) \to \calh^n_H(X;\bfE)$
in the special case $A = \emptyset$.
The functor induced by $\alpha$ on the orbit categories is denoted in the same way
$$\alpha\colon \Or(H) \to \Or(G), \quad H/L \mapsto \ind_{\alpha}(H/L) = G/\alpha(L).$$
There is an obvious natural transformation of functors $\Or(H) \to \GROUPOIDS$
$$
T\colon \calg^H \to \calg^G \circ \alpha.
$$
Its evaluation at $H/L$ is the functor of groupoids
$\calg^H(H/L) \to \calg^G(G/\alpha(L))$ which sends an object $hL$ to the object
$\alpha(h)\alpha(L)$ and a morphism given by $h \in H$ to the morphism $\alpha(h) \in G$.
Notice that $T(H/L)$ is an equivalence if
$\ker(\alpha)$ acts freely on $H/L$. The desired isomorphism 
$$
\ind_{\alpha} \colon H^n_G(\ind_{\alpha} X;\bfE) \to H^n_H(X;\bfE)
$$
is induced by the following map of spectra
\begin{multline*}
\map_{\Or(G)}\left(\map_G(-,\ind_{\alpha}X_+),\bfE \circ \calg^G\right)
\\
\xrightarrow{\cong}
\map_{\Or(G)}\left(\alpha_*(\map_H(-,X_+)),\bfE \circ \calg^G\right)
\\
\xrightarrow{\cong}
\map_{\Or(H)}\left(\map_H(-,X_+),\bfE \circ \calg^G \circ \alpha\right)
\\
\xrightarrow{\map_{\Or(H)}\left(\id,\bfE(T)\right)}
\map_{\Or(H)}\left(\map_H(-,X_+),\bfE \circ \calg^H\right).
\end{multline*}
Here $\alpha_*\map_H(-,X_+)$ is the pointed $\Or(G)$-space which is 
obtained from the pointed $\Or(H)$-space $\map_H(-,X_+)$ by induction, i.e.\
by taking the balanced product over $\Or(H)$ with the 
$\Or(H)$-$\Or(G)$ bimodule $\mor_{\Or(G)}(??,\alpha(?))$
\cite[Definition 1.8]{Davis-Lueck(1998)}.
The second map is given by the adjunction homeomorphism
of induction $\alpha_*$ and restriction $\alpha^*$ (see
\cite[Lemma~1.9]{Davis-Lueck(1998)}). 
The first map comes from the homeomorphism of $\Or(G)$-spaces
$$\alpha_*\map_H(-,X_+) \to \map_G(-,\ind_{\alpha}X_+)$$
which is the adjoint of the obvious map of $\Or(H)$-spaces
$\map_H(-,X_+) \to \alpha^*\map_G(-,\ind_{\alpha}X_+)$ whose evaluation at
$H/L$ is given by $\ind_{\alpha}$.

\em
\end{example}

\typeout{-----------------------  Section 2  ------------------------}

\tit{Modules over a Category}

In this section we give a brief summary about modules over a small
category $\calc$ as far as needed for this paper. They will appear in the definition of
the equivariant Chern character. 

Let $\calc$ be a small category and let $R$ be a
commutative ring.
A \emph{contravariant $R\calc$-module}
is a contravariant  functor from $\calc$ to the category $R\text{ -}\MOD$ of $R$-modules.
Morphisms of contravariant $R\calc$-modules are natural transformations.
Given a group $G$, let $\widehat{G}$ be the category with one object whose  set of
morphisms is given by $G$. Then a contravariant
$R\widehat{G}$-module is the same as a right $RG$-module. Therefore we
can identify the abelian category
$\MOD\text{-}R\widehat{G}$ with the abelian category of right $RG$-modules
$\MOD$-$RG$ in the sequel.
Many of the constructions, which we will introduce for
$R\calc$-modules below,
reduce in the special case $\calc = \widehat{G}$
to their classical versions for $RG$-modules.
The reader should have this example in mind. 
There is  also a covariant version.
In the sequel $R\calc$-module means contravariant $R\calc$-module
unless stated explicitly differently.

The category $\MOD\text{-}R\calc$ of $R\calc$-modules inherits
the structure of an abelian category
from $R\text{ -}\MOD$ in the obvious way, namely
objectwise. For instance a sequence $0 \to M \to N \to P \to 0$ of
contravariant $R\calc$-modules is called \emph{exact} if its evaluation at
each object in $\calc$ is an exact sequence in $R\text{ -}\MOD$. The notion of
an injective and of a  projective $R\calc$-module is now clear.
For a set $S$ denote by $RS$ the free $R$-module with $S$ as basis.
An $R\calc$-module is \emph{free} if it is isomorphic to
$R\calc$-module of the shape $\bigoplus_{i \in I}
R\mor_{\calc}(?,c_i)$ for some index set $I$ and objects $c_i \in \calc$.
Notice that by
the Yoneda-Lemma there is for every $R\calc$-module $N$ and every
object $c$ a bijection of sets
$$\hom_{R\calc}(R\mor_{\calc}(?,c),N) ~ \xrightarrow{\cong} ~ N(c),
\quad \phi ~ \mapsto ~ \phi(\id_x).$$
This implies that every free $R\calc$-module is projective and a
$R\calc$-module is projective if and only if it is a direct summand in a
free $R\calc$-module.
The category of $R\calc$-modules has enough projectives and injectives
(see \cite[Lemma 17.1]{Lueck(1989)} and \cite[Example 2.3.13]{Weibel(1994)}).

Given a contravariant $R\calc$-module $M$ and a covariant $R\calc$-module
$N$, their \emph{tensor product over $R\calc$} is defined to be the following
$R$-module $M \otimes_{R\calc} N$. It is given by
$$M \otimes_{R\calc} N \:=
\bigoplus_{c \in \Ob(\calc)} M(c) \otimes_R N(c)/\sim,$$
where $\sim$ is the typical tensor relation $mf \otimes n = m \otimes fn$,
i.e. for every morphism $f\colon c \to d$ in $\calc$, $m \in M(d)$ and $n \in N(c)$
we introduce the relation $M(f)(m) \otimes n - m \otimes N(f)(n) = 0$.
The main property of this construction is that it is
adjoint to the $\hom_R$-functor in the sense
that for any $R$-module $L$ there are natural isomorphisms of $R$-modules
\begin{eqnarray}
\hom_R(M \otimes_{R\calc} N,L)
& \xrightarrow{\cong} &
\hom_{R\calc}(M,\hom_R(N,L));
\label{adjunction otimes and hom: 1}
\\
\hom_R(M \otimes_{R\calc} N,L)
& \xrightarrow{\cong} &
\hom_{R\calc}(N,\hom_R(M,L)).
\label{adjunction otimes and hom: 2}
\end{eqnarray}

Consider a functor $F\colon \calc \to \cald$. Given a $R\cald$-module $M$,
define its  \emph{restriction with F} to be $F^*M :=M \circ F$.
Given a contravariant $R\calc$-module $M$,
its  \emph{induction with F} is the contravariant 
$R\cald$-module $F_*M$ given by 
\begin{eqnarray}(F_*M)(??) := M(?) \otimes_{R\calc}
  R\mor_{\cald}(??,F(?)),
\label{induction with F}
\end{eqnarray}
and \emph{coinduction with F} is the contravariant 
$R\cald$-module $F_!M$ given by 
\begin{eqnarray}
(F_!M)(??) :=
\hom_{R\calc}(R\mor_{\cald}(F(?),??),M(?)).
\label{coinduction with F}
\end{eqnarray}
Restriction with $F$ can be written  as
$F^*N(?) =  \hom_{R\cald}(R\mor_{\cald}(??,F(?)),N(??))$,
the natural isomorphisms sends $n \in N(F(?))$ to the map
$$R\mor_{\cald}(??,F(?)) \to N(??), \hspace {5mm} 
\phi\colon ?? \to F(?) ~ \mapsto ~ N(\phi)(n).$$
Restriction with $F$ can also be written  as
$F^*N(?) =  R\mor_{\cald}(F(?),??) \otimes_{R\cald} N(??)$,
the natural isomorphisms sends $\phi \otimes_{R\cald} n$ to
$N(\phi)(n)$.
We conclude from \eqref{adjunction otimes and hom: 2}
that $(F_*,F^*)$ and $(F^*,F_!)$ form adjoint pairs, i.e.
for a $R\calc$-module $M$ and a $R\cald$-module $N$ there are natural
isomorphisms of $R$-modules
\begin{eqnarray}
\hom_{R\cald}(F_* M,N) & \xrightarrow{\cong} &
\hom_{R\calc}(M,F^*N); \label{adjunction F_* and F^*}
\\
\hom_{R\cald}(F^*N,M) & \xrightarrow{\cong} &
\hom_{R\calc}(N,F_!M). \label{adjunction F^* and F_!}
\end{eqnarray}

Consider an object $c$ in $\calc$. Let $\aut(c)$ be the group of
automorphism of $c$. We can think of $\widehat{\aut(c)}$ as a
subcategory of $\calc$ in the obvious way. Denote by
$$i(c)\colon \widehat{\aut(c)} \to \calc$$
the inclusion of categories and abbreviate the group ring $R[\aut(c)]$
by $R[c]$ in the sequel. Thus we obtain functors
\begin{eqnarray}
i(c)^* \colon \MOD\text{-}R\calc & \to & \MOD\text{-}R[c];
\label{restriction with i(c)}
\\
i(c)_* \colon  \MOD\text{-}R[c] & \to & \MOD\text{-}R\calc;
\label{induction with i(c)}
\\
i(c)_! \colon  \MOD\text{-}R[c] & \to & \MOD\text{-}R\calc.
\label{coinduction with i(c)}
\end{eqnarray}

The \emph{projective splitting functor}
\begin{eqnarray}
S_c \colon \MOD\text{-}R\calc & \to & \MOD\text{-}R[c]
\label{projective splitting functor S_c}
\end{eqnarray}
sends $M$ to the cokernel of the map
$$\bigoplus_{\substack{f\colon c \to d\\ f \text{ not an isomorphism}}} M(f) \colon
\bigoplus_{\substack{f\colon c \to d\\ f \text{ not an isomorphism}}} M(d) ~ \to ~ M(c).$$
The \emph{injective splitting functor}
\begin{eqnarray}
T_c \colon \MOD\text{-}R\calc & \to & \MOD\text{-}R[c]
\label{injective splitting functor I_c}
\end{eqnarray}
sends $M$ to the kernel of the map
$$\prod_{\substack{f\colon d \to c\\ f \text{ not an isomorphism}}} M(f) \colon
M(c) ~ \to ~ \prod_{\substack{f\colon d \to c\\ f \text{ not an
      isomorphism}}} M(d).$$

From now on suppose that $\calc$ is an EI-category, i.e. a small category such that
endomorphisms are isomorphisms. Then we can define the \emph{inclusion functor}
\begin{eqnarray}
I_c \colon \MOD\text{-}R[c] & \to & \MOD\text{-}R\calc
\label{inclusion functor I_c}
\end{eqnarray}
by $I_c(M)(?) = M \otimes_{R[c]} R\mor(?,c)$ if $c \cong \;?$ in $\calc$ and by $I_c(M)(?)= 0$ otherwise.
Let $B$ be the $R\calc$-$R[c]$-bimodule, covariant over $\calc$ and a right module
over $R[c]$, given by
$$B(c,?) ~ = ~ \begin{array}{lll}
R\mor_{\calc}(c,?) & & \text{ if } c \cong \; ?;
\\
0 & &  \text{ if } c \not\cong \; ?.
\end{array}$$
Let $C$ be the $R[c]$-$R\calc$-bimodule, contravariant over $\calc$ and a left module
over $R[c]$, given by
$$C(?,c) ~ = ~ \begin{array}{lll}
R\mor_{\calc}(?,c) & & \text{ if } c \cong \; ?;
\\
0 & &  \text{ if } c \not\cong \;?.
\end{array}$$
One easily checks that there are natural isomorphisms
\begin{eqnarray*}
S_cM   & \cong & M \otimes_{R\calc} B;
\\
I_c N  & \cong & \hom_{R[c]}(B,N);
\\
T_c M  & \cong & \hom_{R\calc}(C,M);
\\
I_c N  & \cong & N \otimes_{R[c]} C.
\end{eqnarray*}

\begin{lemma}
\label{lem: adjoint and splitting properties}
Let $\calc$ be an EI-category and $c,d$ objects in $\calc$.
\begin{enumerate}

\item \label{lem: adjoint and splitting properties: adjoint}
We obtain adjoint pairs
$(i(c)_*,i(c)^*)$, $(i(c)^*,i(c)_!)$, $(S_c,I_c)$ and $(I_c,T_c)$;

\item \label{lem: adjoint and splitting properties: splitting property}
There are natural equivalences of functors
$S_c \circ i(c)_* \xrightarrow{\cong} \id$ and $T_c \circ i(c)_!  \xrightarrow{\cong} \id$
of functors $\MOD\text{-}R[c] \to \MOD\text{-}R[c]$. If $c \not\cong d$, then
$S_c \circ i(d)_* =  T_c \circ i(d)_! = 0$;

\item \label{lem: adjoint and splitting properties: injectives and projective}
The functors $S_c$ and $i(c)_*$ send projective modules to projective modules.
The functors $I_c$ and $i(c)_!$ send injective modules to injective modules.
\end{enumerate}

\end{lemma}
\begin{proof}
\eqref{lem: adjoint and splitting properties: adjoint} follows from
\eqref{adjunction F_* and F^*}, \eqref{adjunction F^* and F_!} and 
\eqref{adjunction otimes and hom: 1}.
\\[1mm]
\eqref{lem: adjoint and splitting properties: splitting property} This
follows in the case $T_c \circ i(d)_!$ from the following chain of
canonical isomorphisms
\begin{multline*} T_c \circ i(d)_!(M) ~ = ~
\hom_{R\calc}(C(?,c),\hom_{R[d]}(R\mor_{\calc}(d,?),M)) 
\\
\xrightarrow{\cong}
\hom_{R[d]}(C(?,c) \otimes_{R\calc}R\mor_{\calc}(d,?),M) 
\xrightarrow{\cong} \hom_{R[c]}(C(c,d),M),
\end{multline*}
and analogously for $S_d \circ i(c)_*$.
\\[1mm]
\eqref{lem: adjoint and splitting properties: injectives and projective}
The functors $S_c$ and $i(c)_*$ are left adjoint to an exact functor
and hence respect projective. The functors $T_c$ and $i(c)_!$ are right adjoint to an exact functor
and hence respect injective.
\end{proof}

The \emph{length} $l(c) \in \bbN \cup \{\infty\}$ of an object $c$
is the supremum over all natural numbers $l$ for which
there exists a sequence of morphisms
$c_0 \xrightarrow{f_1} c_1 \xrightarrow{f_2} c_2 \xrightarrow{f_3} \ldots
\xrightarrow{f_l} c_l$
such that no $f_i$ is an isomorphism and $c_l = c$.
The \emph{colength} $col(c) \in \bbN \cup \{\infty\}$ of an object $c$
is the supremum over all natural numbers $l$ for which
there exists a sequence of morphisms
$c_0 \xrightarrow{f_1} c_1 \xrightarrow{f_2} c_2 \xrightarrow{f_3} \ldots
\xrightarrow{f_l} c_l$
such that no $f_i$ is an isomorphism and $c_0 = c$.
If each object $c$ has
length $l(c) < \infty$, we say that $\calc$ \emph{has finite length}.
If each object $c$ has
colength $col(c) < \infty$, we say that $\calc$ \emph{has finite colength}.

\begin{theorem}{\bf (Structure theorem for projective and injective $R\calc$-modules).}
\label{the: criterion for projectivity and injectivity}
Let $\calc$ be an EI-category. Then
\begin{enumerate} 
\item \label{the: criterion for projectivity and injectivity:  projectivity} 
Suppose that $\calc$ has finite colength.
Let $M$ be a contravariant $R\calc$-module
such that the $R\aut(c)$-module $S_cM$ is
projective for all objects $c$ in $\calc$.
Let $\sigma_c\colon S_cM \to M(c)$ be an
$R\aut(c)$-section of the canonical projection
$M(c) \to S_cM$. Consider the map of $R\calc$-modules 
\begin{multline*}
\mu(M) \colon \bigoplus_{(c) \in \Is(\calc)} i(c)_*S_cM
\xrightarrow{\bigoplus_{(c) \in \Is(\calc)} i(c)_*\sigma_c}
\bigoplus_{(c) \in \Is(\calc)} i(c)_*M(c)
\\
\xrightarrow{\bigoplus_{(c) \in \Is(\calc)} \alpha(c)} ~ M,
\end{multline*}
where $\alpha(c) \colon i(c)_*M(c) = i(c)_*i(c)^*M \to M$ is the adjoint of the
identity $i(c)^*M \to i(c)^*M$ under the adjunction 
\eqref{adjunction F_* and F^*}.
The map $\mu(M)$ is always surjective. It is bijective if and only if
$M$ is a projective $R\calc$-module;

\item \label{the: criterion for projectivity and injectivity:  injectivity} 
Suppose that $\calc$ has finite length.
Let $M$ be a contravariant $R\calc$-module
such that the $R\aut(c)$-module $T_cM$ is
injective for all objects $c$ in $\calc$.
Let $\rho_c  \colon M(c) \to I_cM $ be an
$R\aut(c)$-retraction of the canonical injection
$T_cM \to M(c)$. Consider the map of $R\calc$-modules 
\begin{multline*}
\nu(M) \colon M \xrightarrow{\prod_{(c) \in \Is(\calc)} \beta(c)} 
\prod_{(c) \in \Is(\calc)} i(c)_!M(c)
\\
\xrightarrow{\prod_{(c) \in \Is(\calc)} i(c)_!\rho_c}
\prod_{(c) \in \Is(\calc)} i(c)_*I_cM
\end{multline*}
where $\beta(c) \colon M \to i(c)_!i(c)^*M = i(c)_!M(c)$ is the adjoint of the
identity $i(c)^*M \to i(c)^*M$ under the adjunction 
\eqref{adjunction F^* and F_!}.
The map $\nu(M)$ is always injective. It is bijective if and only if
$M$ is an injective $R\calc$-module.
\end{enumerate}
\end{theorem}
\begin{proof}
\eqref{the: criterion for projectivity and injectivity:  projectivity}  
A contravariant $R\calc$-module is the same as covariant
$R\calc^{\op}$-module, where $\calc^{\op}$ is the opposite category of $\calc$, just
invert the direction of every morphisms. The corresponding
covariant version of 
assertion~\eqref{the: criterion for projectivity and injectivity:  projectivity}  
is proved in \cite[Theorem~2.11]{Lueck(2002b)}.
\\[2mm]
\eqref{the: criterion for projectivity and injectivity:  injectivity} 
is the dual statement of 
assertion~\eqref{the: criterion for projectivity and injectivity:  projectivity}.
We first show that $\nu(M)$ is always injective. We show by induction over the length
$l(x)$ of an object $x \in \calc$ that $\nu(M)(x)$ is injective.
Let $u$ be an element in the kernel of $\nu(M)(x)$. Consider a morphism
$f \colon y \to x$ which is not an isomorphism. Then $l(y) < l(x)$ and by induction
hypothesis $\nu(M)(y)$ is injective. Since the composite
$\nu(M)(y) \circ M(f)$ factorizes through $\nu(M)(x)$, we have $u \in \ker(M(f))$.
This implies $u \in I_xM$. Consider the composite
$$I_xM \xrightarrow{i} M(x) \xrightarrow{\nu(M)(x)} \prod_{(c) \in \Is(\calc)} i(c)_!I_cM(x)
\xrightarrow{\pr_x} i(x)_!I_xM(x) \xrightarrow{j} I_xM,$$
where $i$ is the inclusion, $\pr_x$ is the projection onto the factor
belonging to the isomorphism class of $x$ and $j$ is the isomorphism
$\hom_{R[x]}(R\mor_{\calc}(x,x),I_xM) \xrightarrow{\cong} I_xM$ sending $\phi$ to
$\phi(\id_x)$. Since this composite is the identity on $I_xM$ and $u$ lies in the kernel of
$\nu(M)(x)$, we conclude $u = 0$.

In particular we see that an injective $R\calc$-module $M$ is trivial if and only if 
$i(d)_!I_dM(x)$ is trivial for all objects $d \in \calc$.

If $\nu(M)$ is bijective and each $I_cM$ is an injective $R[c]$-module,
then $M$ is an injective $R\calc$-module, since $i(c)_!$  sends injective
$R[c]$-modules to injective $R\calc$-modules 
by Lemma~\ref{lem: adjoint and splitting properties}
\eqref{lem: adjoint and splitting properties: injectives and projective}
and the product of injective modules is again injective.

Now suppose that $M$ is injective. Let $N$ be the cokernel of $\nu(M)$.
We have the exact sequence
\begin{eqnarray}
& 0 \to M \xrightarrow{\nu(M)} \prod_{(c) \in \Is(\calc)} i(c)_*I_cM \xrightarrow{\pr} N
\to 0. &
\label{exact sequence of cokernel of nu(M)}
\end{eqnarray}
Since $M$ is injective, this is a split exact sequence of injective
$R\calc$-modules.  Fix an object $d$.
The functors $i(d)_!$ and $I_d$ are left exact and hence send split exact sequences to
split exact sequences. Therefore we obtain a split exact sequence if we apply
$i(d)_!I_d$ to \eqref{exact sequence of cokernel of nu(M)}. Using 
Lemma~\ref{lem: adjoint and splitting properties}
\eqref{lem: adjoint and splitting properties: splitting property}
the resulting exact sequence is isomorphic to the exact sequence
$$0 \to i(d)_!I_dM \xrightarrow{\id} i(d)_!I_dM \to i(d)_!I_dN \to 0.$$
Hence $i(d)_!I_dN$ vanishes for all objects $d$. This implies that $N$ is trivial
and because of \eqref{exact sequence of cokernel of nu(M)} that
$\nu(M)$ is bijective.
\end{proof}

For more details about modules over a category we refer to
\cite[Section 9A]{Lueck(1989)}.

\typeout{-----------------------  Section 3  ------------------------}

\tit{The Associated Bredon Cohomology Theory}

Given a proper equivariant cohomology theory
with values in $R$-modules, we can associate to it another
proper equivariant cohomology theory
with values in $R$-modules satisfying the disjoint union axiom 
called Bredon cohomology,
which is much simpler. The equivariant Chern
character will identify this simpler proper equivariant 
cohomology theory with the given one.

The \emph{orbit category}
$\Or(G)$ has as objects homogeneous spaces $G/H$ and as morphisms $G$-maps.
Let $\Sub(G)$ be the category whose objects are subgroups
$H$ of $G$. For two subgroups $H$ and
$K$ of $G$ denote by $\conhom_G(H,K)$ the set
of group homomorphisms $f\colon H \to K$,
for which there exists an element $g \in G$
with $gHg^{-1} \subseteq K$  such that
$f$ is given by conjugation with $g$, i.e.
$f = c(g)\colon H \to K, \hspace{3mm} h \mapsto ghg^{-1}$.
Notice that $f$ is injective and  $c(g) = c(g')$ holds for two elements $g,g' \in G$ with
$gHg^{-1} \subseteq K$ and $g'H(g')^{-1} \subseteq K$ if and only if
$g^{-1}g'$ lies in the centralizer
$C_GH = \{g \in G \mid gh=hg \mbox{ for all } h \in H\}$
of $H$ in $G$. The group of inner automorphisms of $K$ acts on
$\conhom_G(H,K)$ from the left by composition. Define the set of morphisms
$$\mor_{\Sub(G)}(H,K) ~ := ~ \Inn(K)\backslash \conhom_G(H,K).$$

There is a natural projection $\pr\colon \Or(G) \to \Sub(G)$ which sends
a homogeneous space $G/H$ to $H$.
Given a $G$-map $f\colon G/H \to G/K$, we can choose
an element $g \in G$ with $gHg^{-1} \subseteq K$
and $f(g'H) = g'g^{-1}K$. Then
$\pr(f)$ is represented by $c(g)\colon H \to K$. Notice that
$\mor_{\Sub(G)}(H,K)$ can be identified with the quotient
$\mor_{\Or(G)}(G/H,G/K)/C_GH$,
where $g \in C_GH$ acts on $\mor_{\Or(G)}(G/H,G/K)$
by composition with
$R_{g^{-1}}\colon G/H \to G/H, \hspace{2mm} g'H \mapsto g'g^{-1}H$.

Denote by $\Or(G,\calf) \subseteq \Or(G)$ and
$\Sub(G,\calf) \subseteq \Sub(G)$ the full
subcategories, whose objects $G/H$ and 
$H$ are given by finite subgroups $H \subseteq G$.
Both $\Or(G,\calf)$ and $\Sub(G,\calf)$ are EI-categories of finite length.

Given a proper $G$-cohomology theory $\calh^*_G$ with
values in $R$-modules we obtain for
$n \in \bbZ$ a contravariant $R\Or(G,\calf)$-module
\begin{eqnarray}
 & \calh_G^n(G/?)\colon  \Or(G,\calf) \to R\text{ -}\MOD, \hspace{5mm}
 G/H \mapsto \calh^n_G(G/H). &
\label{Or(G,F)-coefficients H_G(G/?))}
\end{eqnarray}

Let $(X,A)$ be a pair of proper $G$-$CW$-complexes. Then
there is a canonical identification $X^H = \map(G/H,X)^G$.
Thus we obtain contravariant functors
\begin{eqnarray*}
\Or(G,\calf) \to CW\text{-}\PAIRS, & \hspace{5mm} & G/H \mapsto (X^H,A^H);
\\
\Sub(G,\calf) \to CW\text{-}\PAIRS, & \hspace{5mm} &
G/H \mapsto C_GH\backslash(X^H,A^H),
\end{eqnarray*}
where $CW\text{-}\PAIRS$ is the category of pairs of $CW$-complexes. If we compose them
with the covariant functor $CW\text{-}\PAIRS \to \bbZ\text{-}\CHAINCOMPLEXES$
sending $(Z,B)$ to its cellular $\bbZ$-chain complex, then we obtain
the contravariant $\bbZ\Or(G,\calf)$-chain
complex $C^{\Or(G,\calf)}_*(X,A)$ and
the contravariant $\bbZ\Sub(G,\calf)$-chain complex $C^{\Sub(G,\calf)}_*(X,A)$.
Both chain complexes are free in the sense that each chain module
is a free $\bbZ\Or(G,\calf)$-module resp. $\bbZ\Sub(G,\calf)$-module. Namely,
if $X_n$ is obtained from $X_{n-1} \cup A_n$
by attaching the equivariant cells $G/H_i \times D^n$ for $i \in I_n$,
then
\begin{eqnarray}
C^{\Or(G,\calf)}_n(X,A) & \cong & \bigoplus_{i \in I_n}
\bbZ\mor_{\Or(G,\calf)}(G/?,G/H_i);
\label{cellular chain complexes are free: Or}
\\
C^{\Sub(G,\calf)}_n(X,A) & \cong & \bigoplus_{i \in I_n}
\bbZ\mor_{\Sub(G,\calf)}(?,H_i).
\label{cellular chain complexes are free: Sub}
\end{eqnarray}
Given a contravariant $R\Or(G,\calf)$-module $M$,
the \emph{equivariant Bredon cohomology} (see \cite{Bredon(1967a)})
of a pair of proper $G$-$CW$-complexes $(X,A)$ with coefficients
in $M$ is defined by
\begin{eqnarray}
H^n_{\Or(G,\calf)}(X,A;M) & := &
H^n\left(\hom_{\bbZ\Or(G,\calf)}(C^{\Or(G,\calf)}_*(X,A),M)\right).
\label{def. of equi. Bredon cohomology}
\end{eqnarray}
This is indeed a proper $G$-cohomology theory satisfying the disjoint union axiom. 
Hence we can assign to a
proper $G$-homology theory $\calh^*_G$ another proper $G$-cohomology theory
which we call the \emph{associated Bredon cohomology}
\begin{eqnarray}
\calbh^n_G(X,A) & := &
\prod_{p + q = n} H^p_{\Or(G,\calf)}(X,A;\calh^q_G(G/?)).
\label{associated Bredon G-homology theory}
\end{eqnarray}
There is an obvious $\bbZ\Sub(G;\calf)$-chain map
$$\pr_*C_*^{\Or(G,\calf)}(X,A) \xrightarrow{\cong}
C_*^{\Sub(G,\calf)}(X,A)$$
which is bijective because of 
\eqref{cellular chain complexes are free: Or},
\eqref{cellular chain complexes are free: Sub}
and the canonical identification
$$\pr_*\bbZ\mor_{\Or(G,\calf)}(G/?,G/H_i) ~ = ~ \bbZ\mor_{\Sub(G,\calf)}(?,H_i).$$
Given a covariant $\bbZ\Sub(G,\calf)$-module $M$, we get from the
adjunction $ (\pr_*,\pr^*)$ (see Lemma \ref{lem: adjoint and splitting properties}
\eqref{lem: adjoint and splitting properties: adjoint})
natural isomorphisms
\begin{multline}
\hspace{-4mm} H^n_{R\Or(G,\calf)}(X,A;\res_{\pr}M) 
\\ \xrightarrow{\cong}  ~
H^n\left(\hom_{\bbZ\Sub(G,\calf)}\left(C^{\Sub(G,\calf)}_*(X,A),M\right)\right).
\label{ident. of equi. Bredon cohomology}
\end{multline}
This will allow us to work with modules over
the category $\Sub(G;\calf)$ which is smaller
than the orbit category and has nicer properties
from the homological algebra point of view. The main advantage of
$\Sub(G;\calf)$ is that the automorphism groups of every object is finite.

Suppose, we are given a  proper equivariant cohomology theory
$\calh^*_?$ with values in $R$-modules. We get from
\eqref{Or(G,F)-coefficients H_G(G/?))}
for each group $G$ and  $n \in \bbZ$ a covariant $R\Sub(G,\calf)$-module
\begin{eqnarray}
& \calh^n_G(G/?)\colon \Sub(G,\calf) \to R\text{ -}\MOD, \hspace{5mm}
H \mapsto \calh^n_G(G/H).&
\label{Sub(G,F)-coefficients H_G(G/?))}
\end{eqnarray}
We have to show that for $g \in C_GH$ the $G$-map
$R_{g^{-1}}\colon G/H \to G/H, \hspace{3mm} g'H \to g'g^{-1}H$ induces
the identity on $\calh_G^n(G/H)$. This follows from Lemma
\ref{lem: calh_G(G/H) and calh_H(*)}.
We will denote the covariant $R\Or(G,\calf)$-module
obtained by restriction with $\pr\colon \Or(G,\calf) \to \Sub(G,\calf)$ 
from the $R\Sub(G,\calf)$-module
$\calh^n_G(G/?)$ of \eqref{Sub(G,F)-coefficients H_G(G/?))}
again by $\calh^n_G(G/?)$ as introduced already in
\eqref{Or(G,F)-coefficients H_G(G/?))}.

It remains to show that the collection of $G$-cohomology theories
$\calbh^*_G(X,A)$ defined
in \eqref{def. of equi. Bredon cohomology}
inherits the structure of a
proper equivariant cohomology theory, i.e. we have to
specify the induction structure.  We leave it to the reader
to carry out the obvious dualization of the construction for homology in
\cite[Section 3]{Lueck(2002b)} and to check the disjoint union axiom.

\typeout{-----------------------  Section 4  ------------------------}

\tit{The Construction of the Equivariant Cohomological Chern Character}

We begin with explaining the cohomological version of the homological
Chern character due to Dold \cite{Dold(1962)}.

\begin{example}[The non-equivariant Chern character]
 \label{exa: cohomological version of Dold's construction}
\em
Consider a (non-equivariant)
cohomology theory $\calh^*$ with values in $R$-modules. Suppose that $\bbQ \subseteq R$.
For a space $X$ let $X_+$ be the pointed space obtained from $X$ by adding
a disjoint base point $*$. Since the stable homotopy groups
$\pi_p^s(S^0)$ are finite for $p \ge 1$ by results of Serre
\cite{Serre(1953)}, the condition $\bbQ \subseteq R$ imply that the
Hurewicz homomorphism induces isomorphisms
$$\hur_R \colon \pi_p^s(X_+) \otimes_{\bbZ} R 
\xrightarrow{\hur  \otimes_{\bbZ} \id_R} H_p(X) \otimes_{\bbZ} R
\xrightarrow{\cong} H_p(X;R)$$
and that the canonical map
$$\alpha\colon H^p(X;\calh^q(\pt)) \xrightarrow{\cong}
\hom_{\bbQ}(H_p(X;\bbQ),\calh^q(X)) \xrightarrow{\cong}
\hom_R(H_p(X;R),\calh^q(X))$$
is bijective.  Define a map
\begin{eqnarray}
D^{p,q} \colon \calh^{p+q}(X) & \to & \hom_R(\pi^s_p(X_+) \otimes_{\bbZ} R, \calh^q(\pt))
\label{map D^{p,q}}
\end{eqnarray}
as follows. Denote in the sequel by $\sigma^k$ the $k$-fold suspension isomorphism.
Given $a \in \calh^{p+q}(X)$ and an element in  
$\pi^s_p(X_+,*)$ represented by a map $f \colon S^{p+k} \to S^k \wedge X_+$,
we define $D^{p,q}(a)([f]) \in \calh^q(\pt)$ as the image of $a$ under the composite
\begin{multline*}
\calh^{p+q}(X) \xrightarrow{\cong} \widetilde{\calh}^{p+q}(X_+) \xrightarrow{\sigma^k}
\widetilde{\calh}^{p+q+k}(S^k\wedge X_+) \xrightarrow{\widetilde{\calh}^{p+q+k}(f)}
\widetilde{\calh}^{p+q+k}(S^{p+k}) 
\\
\xrightarrow{(\sigma^{p+k})^{-1}}
\widetilde{\calh}^{q}(S^0) \xrightarrow{\cong}  \calh^q(\pt).
\end{multline*}
Then the  (non-equivariant) Chern character for a $CW$-complex $X$ 
is given by the following composite
\begin{multline*}
\ch^n(X)\colon \calh^n(X) \xrightarrow{\prod_{p+q = n} D^{p,q}}
\prod_{p+q = n} 
\hom_R\left(\pi^s_p(X_+,*) \otimes_{\bbZ} R, \calh^q(*)\right) 
\\
\xrightarrow{\prod_{p+q = n} \hom_R(\hur_R^{-1},\id)}
\prod_{p+q = n} \hom_R(H_p(X;R), \calh^q(*)) 
\\
\xrightarrow{\prod_{p+q = n} \alpha^{-1}}
\prod_{p+q = n} H^p(X,\calh^q(*)).
\end{multline*}
There is an obvious version for a pair of $CW$-complexes
$$
\ch^n(X,A) \colon \calh^n(X,A) 
\xrightarrow{\cong}
\prod_{p+q = n} H^p(X,A,\calh^q(*)).
$$
We get a natural transformation $\ch^*$ of cohomology theories with values in $R$-modules.
One easily checks that it is an isomorphism in the case
$X = \pt$.  Hence $\ch^n(X,A)$ is bijective for all relative finite $CW$-pairs $(X,A)$ and
$n \in \bbZ$ by Lemma \ref{lem: role of the disjoint union axiom}
\eqref{lem: role of the disjoint union axiom: transformations of equiv. coho. th.}.
If $\calh^*$ satisfies the disjoint union axiom, then
$\ch^n(X,A)$ is bijective for all $CW$-pairs $(X,A)$ and
$n \in \bbZ$ by Lemma \ref{lem: role of the disjoint union axiom}
\eqref{lem: role of the disjoint union axiom: transformations of equiv. coho. th.}.
\em
\end{example}

Let $R$ be a commutative ring with 
$\bbQ \subseteq R$. Consider an equivariant cohomology theory
$\calh^*_?$ with values in $R$-modules. Let $G$ be a group 
and let $(X,A)$ be a proper $G$-$CW$-pair.
We want to construct an $R$-homomorphism
\begin{multline}
\underline{\ch}^{p,q}_G(X,A)(H)\colon
\calh^{p+q}_G(X,A) 
\\
\to 
\hom_R\left(H_p(C_GH\backslash (X^H,A^H);R),\calh^q_G(G/H)\right)
\label{def of underline{ch}^{p,q}_G}.
\end{multline}
We define it only in the case $A = \emptyset$, the general case is
completely analogous. 
$$\begin{CD}
\calh^{p+q}_G(X)
\\
@V \calh^{p+q}_G(v_H) VV
\\
\calh^{p+q}_G(\ind_{m_H} X^H)
\\
@V \calh^{p+q}_G(\ind_{m_H} \pr_2) VV
\\
\calh^{p+q}_G(\ind_{m_H} EG \times X^H)
\\
@V \ind_{m_H} V\cong V
\\
\calh^{p+q}_{C_GH \times H}( EG \times  X^H)
\\
@V  \left(\ind_{\pr\colon C_GH \times H \to H}\right)^{-1} V \cong V
\\
\calh^{p+q}_H(EG \times_{C_GH} X^H)
\\
@V D^{p,q}_H(EG \times_{C_GH} X^H) VV
\\
\hom_R\left(\pi_p^s((EG \times_{C_GH} X^H)_+) \otimes_{\bbZ} R,\calh^q_H(\pt)\right)
\\
@V \hom_R(\hur_R(EG \times_{C_GH} X^H),\id)^{-1} V V
\\
\hom_R\left(H_p(EG \times_{C_GH} X^H;R),\calh^q_H(\pt)\right)
\\
@ V \hom_R\left(H_p(\pr_1;R),\id\right)^{-1} V V 
\\
\hom_R\left(H_p(C_GH\backslash X^H;R),\calh^q_H(\pt)\right)
\\
@ V \hom_R(\id;(\ind_H^G)^{-1}) V V 
\\
\hom_R\left(H_p(C_GH\backslash X^H;R),\calh^q_G(G/H)\right)
\end{CD}$$
Here are some explanations, more details can be found in
\cite[Section 4]{Lueck(2002b)}. 

We have a left free $C_GH$-action on $EG \times X^H$
by $g(e,x) = (eg^{-1}, gx)$ for $g \in C_GH$, $e \in EG$ and $x \in X^H$. The map
$\pr_1\colon EG \times_{C_GH} X^H \to C_GH\backslash X^H$
is the canonical projection.
Since the projection $BL \to \pt$ induces isomorphisms $H_p(BL;R) \xrightarrow{\cong}
H_p(\pt;R)$ for all $p \in \bbZ$ and finite groups $L$ because of $\bbQ \subseteq R$,  
we obtain for every $p \in \bbZ$ an isomorphism
$$H_p(\pr_1;R)\colon H_p(EG \times_{C_GH} X^H;R) \xrightarrow{\cong}
H_p(C_GH\backslash X^H;R).$$
The group homomorphism $\pr\colon C_GH \times H \to H$
is the obvious projection and
the group homomorphism $m_H\colon C_GH \times H \to G$ sends $(g,h)$ to $gh$.
The $C_GH \times H$-action on $EG \times X^H$
comes from the obvious $C_GH$-action and the trivial $H$-action. 
In particular we equip  $EG \times_{C_GH} X^H$ with the trivial $H$-action. The 
kernels of the two group homomorphisms $\pr$ and $m_H$ act freely on $EG \times X^H$.
We denote by $\pr_2\colon EG \times X^H \to X^H$
the canonical projection.
The $G$-map $v_H\colon \ind_{m_H}X^H =  G \times_{m_H} X^H  \to X$
sends $(g,x)$ to $gx$. 

Since $H$ is a finite group, a $CW$-complex $Z$ equipped with the trivial $H$-action is
a proper $H$-$CW$-complex. Hence we can think of $\calh^*_H$ as an (non-equivariant) homology theory
if we apply it to a $CW$-pair $Z$ with respect to the trivial $H$-action.
Define the map
$$D_H^{p,q}(Z)\colon  \calh_H^{p+q}(Z) \to 
\hom_R(\pi_p^s(Z_+) \otimes_{\bbZ} R,\calh^q_H(\pt))
$$
for a $CW$-complex $Z$ by the map $D^{p,q}$ of \eqref{map D^{p,q}}.

A calculation similar to the one in \cite[Lemma 4.3]{Lueck(2002b)}
shows that the system of maps $\underline{\ch}^{p,q}_G(X,A)(H)$
\eqref{def of underline{ch}^{p,q}_G} fit together to an in $X$ natural
$R$-homomorphism
\begin{multline}
\underline{\underline{\ch}}^{p,q}_G(X,A) \colon \calh^{p+q}_G(X,A)  
\\ \to ~ \hom_{\Sub(G;\calf)}\left(H_p(C_G?\backslash X^?;R),\calh^q_G(G/?)\right).
\label{def of underlineunderline{ch}^{p,q}_G}
\end{multline}

For any contravariant $R\Sub(G;\calf)$-module $M$ and $p \in \bbZ$ there is an in $(X,A)$ natural
$R$-homomorphism
\begin{multline}
\alpha^p_G(X,A;M) \colon H^p_{R\Sub(G;\calf)}(X,A;M)
 ~ \to ~ 
\hom_{\bbQ\Sub(G;\calf)}(H_p(C_G?\backslash X^?;\bbQ),M) 
\\
\xrightarrow{\cong}
\hom_{R\Sub(G;\calf)}(H_p(C_G?\backslash X^?;R),M) 
\label{alpha^n_G(X,A)}
\end{multline}
which is bijective if $M$ is injective as $\bbQ\Sub(G;\calf)$-module.

\begin{theorem}[The equivariant Chern character]
\label{the: equivariant Chern character and injective coeff}
Let $R$ be a commutative ring $R$ with $\bbQ \subseteq R$.
Let $\calh^*_?$ be a proper equivariant cohomology theory with values in $R$-modules.
Suppose that the $R\Sub(G;\calf)$-module
$\calh^q_G(G/?)$ of \eqref{Sub(G,F)-coefficients H_G(G/?))}, 
which sends $G/H$ to $\calh^q_G(G/H)$, is injective as
$\bbQ\Sub(G;\calf)$-module for every group $G$ and every $q \in \bbZ$.

Then we obtain  a  transformation of proper equivariant cohomology theories 
with values in $R$-modules
\begin{eqnarray*}
\ch^*_?\colon \calh^*_? & \xrightarrow{\cong}  & \calbh^*_?,
\end{eqnarray*}
if we define for a group $G$ and a proper $G$-$CW$-pair $(X,A)$
\begin{eqnarray*}
\ch^n_G(X,A) \colon \calh^n_G(X,A)  & \to  &
\calbh^n_G(X,A) :=  \prod_{p+q = n} H^p_{R\Sub(G;\calf)}(X,A;\calh^q_G(G/?))
\end{eqnarray*}
by the composite
\begin{multline*}
\calh^n_G(X,A) \xrightarrow{\prod_{p+q = n}\underline{\underline{\ch}}^{p,q}_G(X,A)} 
\prod_{p+q = n} \hom_{R\Sub(G;\calf)}\left(H_p(C_G?\backslash X^?;R),\calh^q_G(G/?)\right)
\\
\xrightarrow{\prod_{p+q = n} \alpha^p_G(X,A;\calh^q_G(G/?))^{-1}}
\prod_{p+q = n} H^p_{R\Sub(G;\calf)}(X,A;\calh^q_G(G/?))
\end{multline*} 
of the maps defined in \eqref{def of underlineunderline{ch}^{p,q}_G} and
\eqref{alpha^n_G(X,A)}.

The $R$-map $\ch^n_G(X,A)$ is bijective for all proper relative finite $G$-$CW$-pairs $(X,A)$
and $n \in \bbZ$. If $\calh^*_?$ satisfies the disjoint union axiom, then 
the $R$-map $\ch^n_G(X,A)$ is bijective for all proper $G$-$CW$-pairs $(X,A)$
and $n \in \bbZ$. 
\end{theorem}
\begin{proof} 
First one checks that $\ch^*_G$ defines a natural transformation of proper $G$-cohomology
theories. One checks 
for each finite subgroup $H \subseteq G$ and $n \in \bbZ$  that $\ch^n_G(G/H)$ is the
identity if we identify for any $R\Sub(G;\calf)$-module $M$
\begin{multline*}
H^p_{R\Sub(G;\calf)}(G/H;M) =
H^p\left(\hom_{R\Sub(G;\calf)}(C_*^{R\Sub(G;\calf)}(G/H),M\right)
\\
 = ~ \left\{
\begin{array}{lll}
\hom_{R\Sub(G;\calf)}\left(R\mor_{\Sub(G;\calf)}(?,G/H),M\right) = M(G/H) & & \text{ if } p = 0;
\\
0 & & \text{ if } p \not= 0.
\end{array} \right.
\end{multline*}
Finally apply Lemma \ref{lem: role of the disjoint union axiom}
\eqref{lem: role of the disjoint union axiom: transformations of equiv. coho. th.}.
\end{proof}

\begin{remark}[The Atiyah-Hirzebruch spectral sequence for equivariant cohomology]
\label{rem: The Atiyah-Hirzebruch spectral sequence for equivariant cohomology}
\em
There exists a Atiyah-Hirzebruch spectral sequence for equivariant cohomology 
(see \cite[Theroem 4.7 (2)]{Davis-Lueck(1998)}).
It converges to $\calh_G^{p+q}(X,A)$ and has as $E_2$-term the Bredon cohomology groups
$H^p_{R\Sub(G;\calf)}(X,A;\calh^q_G(G/?))$. The conclusion of
Theorem~\ref{the: equivariant Chern character and injective coeff} is that
the spectral sequences collapses.
\em
\end{remark}

\begin{example}[Equivariant Chern character for $\calk^*(G\backslash(X,A))$]
 \label{exa: Chern character for calk^*(G backslas X)} \em
Let $\calk^*$ be a (non-equivariant) cohomology theory 
with values in $R$-modules for a commutative ring with $\bbQ \subseteq R$.
In Example~\ref{exa: cohomology of the quotient or Borel construction}
we have assigned to it an equivariant cohomology theory  by
\begin{eqnarray*}
\calh^n_G(X,A) & = & \calk^n(G\backslash (X,A)).
\end{eqnarray*}
We claim that the assumptions appearing in 
Theorem~\ref{the: equivariant Chern character and injective coeff} are satisfied
We have to show that the constant functor
$$\underline{\calk^q(\pt)} \colon \Sub(G;\calf) \to \bbQ\text{ -}\MOD, 
\quad H ~ \mapsto \calk^q(\pt)$$
is injective. Let $i \colon \Sub(\{1\}) \to \Sub(G;\calf)$ the obvious inclusion of
categories. Since the object $\{1\}$ is an initial object in
$\Sub(G;\calf)$, the $\bbQ\Sub(G;\calf)$-modules $i_!(\calk^q(\pt))$
and $\underline{\calk^q(\pt)}$ are isomorphic.
Since $i_!$ sends an injective $\bbQ$-module to an injective $R\Sub(G;\calf)$-module
by Lemma\ref{lem: adjoint and splitting properties}
\eqref{lem: adjoint and splitting properties: injectives and projective}
and $\calh^q(\pt)$ is injective as $\bbQ$-module,
$\underline{\calk^q(\pt)}$ is injective as $\bbQ\Sub(G;\calf)$-module. From 
Theorem~\ref{the: equivariant Chern character and injective coeff} we get a
transformation of equivariant cohomology theories
\begin{multline*}
\ch^n_G(X,A) \colon \calk^n(G\backslash (X,A))   ~ \xrightarrow{\cong} ~
\prod_{p+q = n} H^p_{R\Sub(G;\calf)}(X,A;\underline{\calh^q(\pt)}) 
\\ ~ = ~ \prod_{p+q = n} H^p(G\backslash (X,A);\calh^q(\pt)).
\end{multline*}
One easily checks that this is precisely the Chern character of 
Example~\ref{exa: cohomological version of Dold's construction} applied to
$\calk^*$ and the $CW$-pair $G\backslash(X,A)$. \em
\end{example}

\typeout{-----------------------  Section 5  ------------------------}

\tit{Mackey Functors}

In Theorem~\ref{the: equivariant Chern character and injective coeff}
the assumption appears that the contravariant $R\Sub(G;\calf)$-module
$\calh^q_G(G/?)$ is injective for each $q \in \bbZ$. We want to give a
criterion which ensures that this assumption is satisfies and which turns out
to apply to all cases of interest.

Let $R$ be a commutative ring.
Let $\FGINJ$ be the category of finite
groups with injective group homomorphisms as morphisms.
Let $M\colon \FGINJ \to R\text{ -}\MOD$ be a bifunctor, i.e. a pair $(M_*,M^*)$
consisting of a covariant functor $M_*$ and a contravariant
functor $M^*$ from $\FGINJ$ to $R\text{ -}\MOD$ which agree on objects.
We will often denote for an injective  group homomorphism $f\colon H \to G$
the map $M_*(f)\colon M(H) \to M(G)$ by $\ind_f$ and
the map $M^*(f)\colon M(G) \to M(H)$ by $\res_f$ and write
$\ind_H^G = \ind_f$ and $\res_G^H = \res_f$ if $f$ is an inclusion of
groups. We call such a bifunctor $M$ a \emph{Mackey functor}
with values in $R$-modules if
\begin{enumerate}
\item
For an inner automorphism $c(g)\colon G \to G$ we have
$M_*(c(g)) = \id\colon M(G) \to M(G)$;

\item
For an isomorphism of groups $f\colon G \xrightarrow{\cong} H$
the composites $\res_f \circ \ind_f$ and $\ind_f \circ \res_f$
are the identity;

\item Double coset formula\\[1mm]
We have for two subgroups $H,K \subseteq G$
$$\res_G^K \circ \ind_H^G ~ = ~\sum_{KgH \in K\backslash G/H}
\ind_{c(g)\colon H\cap g^{-1}Kg \to K}
\circ \res_{H}^{H\cap g^{-1}Kg},$$
where $c(g)$ is conjugation with $g$, i.e. $c(g)(h) = ghg^{-1}$.
\end{enumerate}

Let $G$ be a group. In the sequel we denote for a subgroup 
$H \subseteq G$ by $N_GH$
the normalizer and by $C_GH$ the centralizer of $H$ in $G$ and
by $W_GH$ the quotient $N_GH/H\cdot C_GH$. Notice that
$W_GH$ is finite if $H$ is finite.
Let $R$ be a commutative ring. Let $M$ be a Mackey functor with values in $R$-modules.
It induces a contravariant $R\Sub(G,\calf)$-module denoted in the same way
$$M\colon  \Sub(G,\calf) \to R\text{ -}\MOD, \hspace{5mm} (f\colon H \to K) \mapsto
\left(M^*(f) \colon M(H) \to M(K)\right).$$
We want to use Theorem~\ref{the: criterion for projectivity and injectivity} 
\eqref{the: criterion for projectivity and injectivity:  injectivity} 
to show that $M$ is injective and analyse its structure.
The $R[W_GH]$-module $T_HM$ introduced in 
\eqref{injective splitting functor I_c} is the same
as the kernel of 
$$\prod_{K \subsetneq H}
 ~ M(i_K) \colon  M(H) ~ \to ~ 
\prod_{K \subsetneq H} ~M(K),$$
where for each subgroup $K \subsetneq H$ different from $H$ 
we denote by $i_K$ the inclusion.
Suppose that $R[W_GH]$-module $T_HM$ is injective for every finite
subgroup $H \subseteq G$.  For every finite
subgroup $H \subseteq G$ choose a retraction $\rho_H\colon M(H) \to T_HM$ of the
inclusion $T_HM \to M(H) $. Denote by $I= \Is(\Sub(G,\calf))$ the set of isomorphism
classes of objects in $\Sub(G;\calf)$ which is the same as the set 
of conjugacy classes $(H)$ of finite subgroups $H$ of $G$.
Let 
\begin{eqnarray}
\nu =\nu(M) \colon M & \to &   \prod_{(K) \in I} i(K)_! \circ T_K(M)
\label{def of map nu}
\end{eqnarray}
be the homomorphism of $R\Sub(G,\calf)$-modules 
uniquely determined by the property that for any $(K) \in I$
its composition with the projection onto the factor indexed by $(K)$
is the adjoint of $\rho_K\colon M(K) \to T_KM$ for the adjoint pair $(i(K)^*,i(K)_!)$.

\begin{theorem}[Injectivity and Mackey functors] \label{the: injectivity of a Mackey functor}
Let $G$ be a group and let $R$ be a commutative ring such that the order
of every finite subgroup of $G$ is invertible in $R$. 
Suppose that the $R[W_GH]$-module $T_HM$ is injective for each finite
subgroup $H \subseteq G$. Then
$M$ is injective as $R\Sub(G,\calf)$-module and the map $\nu$ of
\eqref{def of map nu} is bijective.
\end{theorem}
\begin{proof} 
The map $\nu$ of \eqref{def of map nu} is the map $\nu(M)$ appearing in 
Theorem~\ref{the: criterion for projectivity and injectivity} 
\eqref{the: criterion for projectivity and injectivity:  injectivity}. 
Because of Theorem~\ref{the: criterion for projectivity and injectivity} 
\eqref{the: criterion for projectivity and injectivity:  injectivity} it
suffices to show for each finite subgroup $H \subseteq G$ 
that $\nu(M)(H)$ is surjective.

Fix for any $(K) \in I$ a representative $K$. Then choose for any
$W_GK \cdot f\in W_GK\backslash\mor(K,H)$ an element $f \in
\conhom(K,H)$ which represents a morphism $f \colon K \to H$ in
$\Sub(G;\calf)$ which belongs to $W_GK \cdot f\in W_GK\backslash\mor(K,H)$.
Notice that $W_GK$ is the automorphism group of the object $K$ in
$\Sub(G;\calf)$ and $W_GK$, $\mor(K,H)$ and 
$W_GK\backslash\mor(K,H)$ are finite. 
With these choices we get for every object $H$ in $\Sub(G;\calf)$ an identification
$$i(K)_!T_KM(H) ~ = ~ \hom_{RW_GK}(R\mor(K,H),T_KM) = 
\prod_{\substack{W_GK \cdot f \in \\ W_GK\backslash\mor(K,H)}} T_KM^{W_GK_f}$$
where $W_GK_f \subseteq W_GK$ is the isotropy group of $f$ under the
$W_GK$-action on $\mor(K,H)$.
Under this identification $\nu(H)$ becomes the map
$$\nu(H) \colon M(H) \to \prod_{(K) \in I} ~  
\prod_{\substack{W_GK \cdot f\in \\  W_GK\backslash\mor(K,H)}} T_KM^{W_GK_f}$$
for which the component of $\nu(H)(m)$, which belongs $(K) \in I$ and
$W_GK \cdot f\in W_GK\backslash\mor(K,H)$, is $\rho_K \circ \res_f(m)$ for $m \in M(H)$. 
Notice that the image of $\res_f$ always is contained in $M(K)^{W_GK_f}$.
Next we define a map
$$\mu(H)\colon \bigoplus_{\substack{(K) \in I,\\ (K) \le (H)}}  ~ \bigoplus_{\substack{W_GK \cdot
f\in \\ W_GK\backslash\mor(K,H)}} T_KM^{W_GK_f} \to M(H)$$
by requiring that its restriction to the summand, which belongs to $(K) \in I$ and 
$W_GK \cdot f\in  W_GK\backslash\mor(K,H)$, is the composite of the inclusion
$T_KM^{W_GK_f} \to M(K)$ with $\ind_f\colon M(K) \to M(H)$.
We want to show that the composite
\begin{multline*}
\nu(H) \circ \mu(H) \colon \bigoplus_{\substack{(K) \in I,\\ (K) \le (H)}} ~
 \bigoplus_{\substack{W_GK \cdot f\in\\ W_GK\backslash\mor(K,H)}} T_KM^{W_GK_f} \\
\to 
\prod_{(K) \in I}  ~ 
\prod_{\substack{W_GK \cdot f\in\\ W_GK\backslash\mor(K,H)}} T_KM^{W_GK_f} 
~ = ~ 
\bigoplus_{\substack{(K) \in I,\\ (K) \le (H)}}  ~
\bigoplus_{\substack{W_GK \cdot f\in\\W_GK\backslash\mor(K,H)}} T_KM^{W_GK_f}
\end{multline*}
is bijective. If $K$ is subconjugated to $H$, we write $(K) \le (H)$. 
Fix $(K), (L) \in I$ with $(K) \le (H)$ and $(L) \le (H)$ and 
$W_GK\cdot f \in \mor(K,H)$ and $W_GL\cdot g \in \mor(L,H)$.  Then 
the homomorphism $T_{K}M^{{W_GK}_{f}} \to T_{L}M^{{W_GL}_{g}}$ given by $\nu(H) \circ \mu(H)$ 
and the summands corresponding to $(K,f)$ and
$(L, g)$ is induced by the composite
\begin{multline}
\alpha_{(K,f),(L,g)} \colon T_{K}M^{{W_GK}_{f}} \xrightarrow{i} M(K) 
\xrightarrow{\ind_{f \colon K \to \im(f)}} 
M(\im(f)) \xrightarrow{\ind_{\im(f)}^H} M(H) 
\\
\xrightarrow{\res_H^{\im(g)}} M(\im(g)) 
\xrightarrow{\res_{g \colon L \to \im(g)}} M(L)
\xrightarrow{\rho_{L}} T_{L}M,
\label{alpha_{(K,f),(L,g)}} 
\end{multline}
where $i$ is the inclusion. The double coset formula implies
\begin{multline}
\res_H^{\im(g)} \circ \ind_{\im(f)}^H 
\\
 = ~ 
\sum_{\substack{\im(g)h\im(f) \in \\ \im(g)\backslash H/\im(f)}}
\ind_{c(h)\colon \im(f)\cap h^{-1}\im(g)h \to \im(g)}
 \circ \res_{\im(f)}^{\im(f)\cap h^{-1}\im(g)h}.
\label{formula for res_H^{im(g)} circ ind_{im(f)}^H}
\end{multline}
The composite 
\begin{multline*}
T_{K}M^{{W_GK}_{f}} \xrightarrow{i} M(K) 
\xrightarrow{\ind_{f \colon K \to \im(f)}} 
M(\im(f)) 
\\
\xrightarrow{\res_{\im(f)}^{\im(f)\cap h^{-1}\im(g)h}} 
M(\im(f)\cap h^{-1}\im(g)h)
\end{multline*}
is trivial by the definition of $T_{K}M$ if
$\im(f)\cap h^{-1}\im(g)h \not= \im(f)$ holds. 
Hence $\alpha_{(K,f),(L,g)} \not=0$ is only possible if
$\im(f)\cap h^{-1}\im(g)h = \im(f)$ for some $h \in H$ and hence 
$(K) \le  (L)$ hold.

Suppose that $(K) = (L)$. Then $K = L$ by our choice of representatives. Suppose
that $\alpha_{(K,f),(K,g)} \not=0$. We have already seen that
this implies $\im(f)\cap h^{-1}\im(g)h = \im(f)$ for some $h \in
H$. Since $|\im(f)| = |K| = |\im(g)|$ we conclude 
$h^{-1}\im(g)h = \im(f)$ and therefore $W_GK \cdot f = W_GK\cdot g$ in 
$W_GK\backslash\mor(K,H)$. This implies already $f = g$ as group
homomorphism $K \to H$ by our choice of representatives.
The double coset formula~\eqref{formula for res_H^{im(g)} circ ind_{im(f)}^H} 
implies that $\alpha_{(K,f),(K,f)}$ is 
$|H \cap N_G\im(f)| \cdot \id_{T_{K}M^{{W_GK}_{f}}}$ since for all 
$h \in N_G\im(f) \cap H$ the composite
\begin{multline*}
T_{K}M^{{W_GK}_{f}} \xrightarrow{i} M(K) 
\xrightarrow{\ind_{f \colon K \to \im(f)}} 
M(\im(f)) 
\\
\xrightarrow{\ind_{c(h)\colon \im(f) \to \im(f)}}
M(\im(f))
\end{multline*}
agrees with 
$T_{K}M^{{W_GK}_{f}} \xrightarrow{i} M(K) 
\xrightarrow{\ind_{f \colon K \to \im(f)}} M(\im(f))$. Since the order of
$|H \cap N_G\im(f)|$ is invertible in $R$ by assumption,
$\alpha_{(K,f),(K,f)}$ is bijective.  

We conclude that $\nu(H) \circ \mu(H)$ can be written as a matrix of maps
which has upper triangular form and isomorphisms on the diagonal.
Therefore $\nu(H) \circ \mu(H)$ is surjective. This shows that $\nu(H)$ is surjective.
This finishes the proof of Theorem \ref{the: injectivity of a Mackey functor}. 
\end{proof}

\begin{theorem}[The equivariant Chern character and Mackey structures]
\label{the: Chern for calh^*_G with a Mackey structure}
Let $\calh^*_?$ be a proper equivariant cohomology theory. 
Define a contravariant functor
$$\calh^q_?(\pt) \colon \FGINJ \to R\text{ -}\MOD$$
by sending a homomorphism $\alpha \colon H \to K$ to the composite 
$$\calh_K^q(\pt) \xrightarrow{\calh^q(\pr)}
\calh^q_K(K/H) \xrightarrow{\ind_{\alpha}} \calh^q_H(\pt)$$ 
where $\pr \colon H/K = \ind_{\alpha}(\pt) \to \pt$ is
the projection and $\ind_{\alpha}$ comes from the induction 
structure of $\calh^*_?$. Suppose that it extends to a Mackey functor
for every $q \in \bbZ$. Then

\begin{enumerate}

\item \label{the: Chern for calh^*_G with a Mackey structure: injectivity}
For every group $G$ the $R\Sub(G;\calf)$-module $\calh^q_G(G/?)$  of 
\eqref{Sub(G,F)-coefficients H_G(G/?))} is injective as $R\Sub(G;\calf)$-module,
provided that $R$ is semisimple;

\item \label{the: Chern for calh^*_G with a Mackey structure: calh cong calbh}
We obtain a natural transformation of proper
equivariant cohomology theories with values in $R$-modules
$$
\ch^*_?(X,A) \colon \calh^*_?  ~ \to ~ \calbh^*_?.$$
In particular we get for every  proper $G$-$CW$-pair $(X,A)$ 
and every $n \in \bbZ$ a natural $R$-homomorphism
\begin{multline*}
\ch^n_G(X,A) \colon \calh^n_G(X,A)   
\\
~ \to ~
\calbh^n_G(X,A) ~ := ~ \prod_{p+q = n} H^p_{R\Sub(G;\calf)}(X,A;\calh^q_G(G/?)).
\end{multline*}
It is bijective for all proper relative finite $G$-$CW$-pairs $(X,A)$ and $n \in \bbZ$.
If $\calh^*_?$ satisfies the disjoint union axiom, 
it is bijective for all proper $G$-$CW$-pairs $(X,A)$ and $n \in \bbZ$;

\item \label{the: Chern for calh^*_G with a Mackey structure: splitting}
Define for finite subgroup $H \subseteq G$ the $R[W_GH]$-module
$T_H(\calh^q_H(\pt))$ by 
$$
 \ker\left(\prod_{L \subsetneq H} 
\ind_L^K \circ \calh^q(\pr \colon H/L \to \pt) \colon \calh^q_H(\pt) \to
  \prod_{L \subsetneq H} \calh^q_L(\pt)\right).
$$
Then the Bredon cohomology
$\calbh^n_G(X,A)$ of a proper $G$-$CW$-pair $(X,A)$
is naturally $R$-isomorphic to 
\begin{multline*}
\prod_{p +q = n} ~ \prod_{(H), H \subseteq G \text{ finite}}~ 
\hom_{RW_GH}\left(H_p(C_GH\backslash X^H;R),T_H(\calh^q_H(\pt)\right)
\end{multline*}
\end{enumerate}
provided that $R$ is semisimple.
\end{theorem}
\begin{proof}
\eqref{the: Chern for calh^*_G with a Mackey structure: injectivity} This follows from
Theorem \ref{the: injectivity of a Mackey functor} since  for every finite subgroup $H
\subseteq G$ the group $W_GH$ is finite and hence the ring $R[W_GH]$ is semisimple and 
every $R[W_GH]$-module is injective.
\\[2mm]
\eqref{the: Chern for calh^*_G with a Mackey structure: calh cong calbh} This follows from
assertion \eqref{the: Chern for calh^*_G with a Mackey structure: injectivity}
applied in the case $R = \bbQ$ 
together with Theorem \ref{the: equivariant Chern character and injective coeff}.
\\[2mm]
\eqref{the: Chern for calh^*_G with a Mackey structure: splitting} 
Since $R$ is semisimple,  the ring $R[W_GH]$ is semisimple and 
every $R[W_GH]$-module is injective for every finite subgroup $H$.
Because the map $\nu$ of \eqref{def of map nu} is an isomorphism by 
Theorem \ref{the: injectivity of a Mackey functor}, it remains to show
for a $C_GH$-module $N$
$$\hom_{R\Sub(G;\calf)}(H_p(C_G?\backslash X^?;R),i(H)_!N) ~ = ~  
\hom_{RW_GH}(H_p(C_GH\backslash X^H;R),N).$$
This follows from the adjunction $(i(H)^*,i(H)_!)$  of
Lemma \ref{lem: adjoint and splitting properties}
\eqref{lem: adjoint and splitting properties: adjoint}. 
\end{proof}

\begin{example}[Mackey structures for Borel cohomology]
\label{exa: Mackey structures for Borel cohomology}
\em
Let $\calk^*$ be a cohomology theory for (non-equivariant)
$CW$-pairs with values in $R$-modules for a commutative ring
$R$ such that $\bbQ \subseteq R$ and $R$ is semisimple.
In Example~\ref{exa: cohomology of the quotient or Borel construction}
we have assigned to it an equivariant cohomology theory called equivariant Borel
cohomology by
\begin{eqnarray*}
\calh^n_G(X,A) & = & \calk^n(EG \times_G (X,A)).
\end{eqnarray*}
We claim that the assumptions appearing in 
Theorem~\ref{the: Chern for calh^*_G with a Mackey structure} are satisfied.
Namely, the contravariant functor
$$\FGINJ \to R\text{ -}\MOD, \quad H ~ \mapsto ~ \calk^n(BH)$$
extends to a Mackey functor, the necessary covariant functor comes from the
Becker-Gottlieb transfer (see for instance \cite{Feshbach(1979)} and
\cite[Corollary 6.4 on page 206]{Lewis-May-Steinberger(1986)}).
Hence we get from Theorem~\ref{the: Chern for calh^*_G with a Mackey structure}
for every group $G$ and every proper $G$-$CW$-pair $(X,A)$ natural $R$-maps
\begin{multline*}
\ch^n_G(X,A) \colon \calk^n(EG \times_G(X,A))   
~ \xrightarrow{\cong}  ~
\prod_{p+q = n} H^p_{R\Sub(G;\calf)}(X,A;\calk^q(B?))
\\
\cong 
\prod_{p +q = n} ~ \prod_{(H), H \subseteq G \text{ finite}}~ 
\hom_{RW_GH}(H_p(C_GH\backslash X^H;R),T_H(\calk^q(BH)),
\end{multline*}
if we define 
\begin{eqnarray*}
T_H(\calk^q(BH)) ~ := ~
\ker\left(\prod_{K \subsetneq H} \calk^q(BK \to BH) \colon \calk^q(BH) \to
  \prod_{K \subsetneq H} \calk^q(BK)\right).
\end{eqnarray*}
If $(X,A)$ is relative finite or if $\calk^*$ satisfies the disjoint union axiom, then these
maps $\ch^n_G(X,A)$ are bijective.
\end{example}

\begin{remark} \label{Theorem 0.1 does not follow} \em
We mention that this does \emph{not} prove 
Theorem~\ref{the: rational computation of K^*(BG)} since
we cannot apply it to $\calk^* := K^* \otimes_{\bbZ} \bbQ$. The problem is
that $K^* \otimes_{\bbZ} \bbQ$ defines all axioms of a cohomology theory but not the
disjoint union axiom. But this is needed if we want to deal
with classifying spaces $BG$ of groups which are not finite $CW$-complexes,
for instance of groups containing torsion (see also 
Remark~\ref{rem: The disjoint union axiom is not compatible with - otimes_{bbZ} bbQ}
and Example~\ref{exa: Rationalizing topological  K -theory}).

A proof of Theorem~\ref{the: rational computation of K^*(BG)} will be given
in \cite{Lueck(2004i)}.
\em
\end{remark}

\begin{example}[Equivariant $K$-theory and Mackey structures] 
\label{exa: Equivariant K-theory and Mackey structures} \em
In Example~\ref{exa: equivariant K-theory} we have introduced the equivariant
cohomology theory $K^*_?$ given by topological $K$-theory. Recall that it takes values
in $R$-modules for $R = \bbZ$. Notice that for a finite group $H$ we get an identification
of $K^0_H(\pt)$ with the complex representation ring $R(H)$ and the associated
contravariant functor
$$K^q_? \colon \FGINJ \to R\text{ -}\MOD, \quad H ~ \mapsto ~ K^q_H(\pt) = R(?)$$
sends an injective group homomorphism $\alpha \colon H \to G$ of finite groups
to the homomorphism of abelian groups $R(G) \to R(H)$ given by restriction with $\alpha$.
Induction with $\alpha$ induces a covariant functor $H \mapsto R(H)$ and it turns
out that this defines a Mackey structure on $K^q_?$. 

For rationalized equivariant topological $K$-theory $K^*_?\otimes_{\bbZ} \bbQ $ 
the equivariant Chern character of
Theorem~\ref{the: Chern for calh^*_G with a Mackey structure} can be identified 
with the one constructed in \cite{Lueck-Oliver(2001b)}
for proper relative finite $G$-$CW$-pairs $(X,A)$.
\em
\end{example}

\typeout{-----------------------  Section 6  ------------------------}

\tit{Multiplicative Structures}

Next we want to introduce a multiplicative structure on a proper equivariant
cohomology theory $\calh^*_?$ and show that it induces one on
the associated Bredon cohomology $\calbh^*_?$ such that the
equivariant Chern character is compatible with it.

We begin with the non-equivariant case. Let $\calh^*$ be
a (non-equivariant) cohomology theory with values in $R$-modules. 
A \emph{multiplicative structure} assigns to a $CW$-complex $X$ with
$CW$-subcomplexes $A,B \subseteq X$ natural $R$-homomorphisms
\begin{eqnarray}
\cup \colon \calh^{n}(X,A) \otimes_R \calh^{n'}(X,B) & \to &
\calh^{n+n'}(X,A\cup B).
\label{non-equivariant cup -product}
\end{eqnarray}
This product is required to be compatible with the boundary
homomorphism of  the long exact sequence
of a pair, to be graded commutative, to be 
associative  and to have a unit $1 \in \calh^0(\pt)$.

Given a multiplicative structure on $\calh^*$, we obtain for every
$p,q \in \bbZ$ a pairing
$$\cup \colon \calh^q(\pt) \otimes_R \calh^{q'}(\pt) ~ \to ~ \calh^{q+q'}(\pt).$$
It yields on singular (or equivalently cellular) cohomology a product
$$H^p(X,A;\calh^q(\pt)) \otimes_R H^{p'}(X,B;\calh^{q'}(\pt)) ~ \to 
H^{p+p'}(X,A\cup B;\calh^{q+q'}(\pt)).$$
The collection of these pairings induce a multiplicative structure
on the cohomology theory given by 
$\prod_{p+q = n} H^p(X,A;\calh^q(\pt))$.  The straightforward proof of
the next lemma is left to the reader.

\begin{lemma} \label{lem: Dold's construction is comp with cup}
Let $R$ be a commutative ring with $\bbQ \subseteq R$. 
Let $\calh^*$ be a (non-equivariant) cohomology theory satisfying
the disjoint union axiom which comes with a multiplicative structure. 

Then the (non-equivariant) Chern character of 
Example~\ref{exa: cohomological version of Dold's construction} 
$$
\ch^n(X,A)\colon \calh^n(X,A) 
\xrightarrow{\cong}
\prod_{p+q = n} H^p(X,A,\calh^q(*))
$$
is compatible with the given multiplicative structure on $\calh^*$ and
the induced multiplicative structure on the target.
\end{lemma}

Next we deal with the equivariant version. We only deal with the
proper case, the definitions below make also sense without this condition.

Let  $\calh^*_G$ be a proper $G$-cohomology theory. A \emph{multiplicative structure} 
assigns to a  proper $G$-$CW$-complex $X$ with
$G$-$CW$-subcomplexes $A,B \subseteq X$ natural $R$-homomorphisms
\begin{eqnarray}
\cup \colon \calh^{n}_G(X,A) \otimes_R \calh^{n'}_G(X,B) & \to &
\calh^{n+n'}_G(X,A\cup B).
\label{equivariant cup -product}
\end{eqnarray}
This product is required to be compatible with the boundary
homomorphism of  the long exact sequence of a $G$-$CW$-pair, to be graded commutative,
to be associative and to have a unit $1 \in \calh^0_G(X)$ for every proper
$G$-$CW$-complex $X$

Let $\calh^*_?$ be a proper equivariant cohomology theory. A
\emph{multiplicative structure} on it assigns a multiplicative
structure to the associated proper $G$-coho\-mo\-logy theory $\calh^*_G$ for
every group $G$ such that for each group homomorphism
$\alpha \colon H \to G$ the maps given by the induction structure of \eqref{induction structure}
\begin{eqnarray*}
\ind_{\alpha}\colon \calh^n_G(\ind_{\alpha}(X,A)) 
&\xrightarrow{\cong} &
\calh^n_H(X,A)
\end{eqnarray*}
are in the obvious way compatible with the multiplicative structures on $\calh^*_G$ and
$\calh^*_H$. 

Next we explain how a given multiplicative structure on $\calh^*_?$
induces one on $\calbh^*_?$. We have to specify for every group $G$
a multiplicative structure on the $G$-cohomology theory
$\calbh^*_G$. Consider a $G$-$CW$-complex $X$ with
$G$-$CW$-subcomplexes $A,B \subseteq X$. For two contravariant $R\Or(G;\calf)$-chain
complexes $C_*$ and $D_*$ define the contravariant $R\Or(G;\calf)$-chain
complexes $C_* \otimes_R D_*$ by sending $G/H$ to the tensor product
of $R$-chain complexes $C_*(G/H) \otimes_R D_*(G/H)$. Let
$$a_* \colon C^{R\Or(G;\calf)}_*(X,A) \otimes_R C^{R\Or(G;\calf)}_*(X,B) ~
\xrightarrow{\cong}  C^{R\Or(G;\calf)}_*((X,A) \times (X,B))$$
be the isomorphism of $R\Or(G;\calf)$-chain complexes which is given
for an object $G/H$ by the natural isomorphism of cellular $R$-chain
complexes 
$$C_*(X^H,A^H) \otimes_R C_*(X^H,B^H) \xrightarrow{\cong}
C_*((X^H,A^H) \times (X^H,B^H)).$$
The multiplicative structure on $\calh_G^*$ yields a 
map of contravariant $R\Or(G;\calf)$-modules
$$c \colon \calh^q_G(G/?) \otimes_R \calh^q_G(G/?) \to \calh^{q+q}_G(G/?).$$
Let
$$\Delta \colon (X;A\cup B) ~ \to (X,A) \times (X,B), \quad x \mapsto
(x,x)$$
be the diagonal embedding. Define a $R$-cochain map by the composite
\begin{multline*}
b^* \colon
\hom_{R\Or(G;\calf)}\left(C^{R\Or(G;\calf)}_*(X,A),\calh^q_G(G/?)\right)
\\
\otimes_R \hom_{R\Or(G;\calf)}\left(C^{R\Or(G;\calf)}_*(X,B),\calh^{q'}_G(G/?)\right)
~ \xrightarrow{\otimes_R} 
\\
\hom_{R\Or(G;\calf)}\left(C^{R\Or(G;\calf)}_*(X,A)\otimes_R C^{R\Or(G;\calf)}_*(X,B),
\calh^q_G(G/?) \otimes_R \calh^{q'}_G(G/?)\right)
\\
\xrightarrow{\hom_{R\Or(G;\calf)}\left((a_*)^{-1},c\right)}
\hom_{R\Or(G;\calf)}\left(C^{R\Or(G;\calf)}_*((X,A)\times (X,B)),
  \calh^{q+q'}_G(G/?)\right)
\\
\xrightarrow{\hom_{R\Or(G;\calf)}(C^{R\Or(G;\calf)}_*(\Delta),\id)}
\hom_{R\Or(G;\calf)}\left(C^{R\Or(G;\calf)}_*(X,A \cup B),\calh^{q+q'}_G(G/?)\right).
\end{multline*}
There is a canonical $R$-map
$$H^*(C^* \otimes_R D^*) \to H^*(C^* \otimes_R D^*)$$
for two $R$-cochain complexes $C^*$ and $D^*$. This map together with
the map induced by $b^*$ on cohomology yields an $R$-homomorphism
\begin{multline*}
H^p_{R\Or(G;\calf)}\left(X,A;\calh^q_G(G/?)\right) \otimes_R 
H^{p'}_{R\Or(G;\calf)}(X,B;\calh^{q'}_G(G/?)) 
\\  \to ~ 
H^{p+p'}_{R\Or(G;\calf)}\left(X,A\cup B;\calh^q_G(G/?)\right).
\end{multline*}
The collection of these $R$-homomorphisms yields 
the desired multiplicative structure on $\calbh^*_G$.
We leave it to the reader to check that the axioms of a multiplicative
structure on $\calbh^*_G$ are satisfied and that all these are
compatible with the induction structure so that we obtain a
multiplicative structure on the equivariant cohomology theory $\calbh^*_?$.

We also omit the lengthy but straightforward proof of the following
result which is based on 
Theorem~\ref{the: equivariant Chern character and injective coeff},
Lemma ~\ref{lem: Dold's construction is comp with cup} and the
compatibility of the multiplicative structure with the induction
structure.

\begin{theorem}[The equivariant Chern character and multiplicative structures]
\label{the: equivariant Chern character and multiplicative structures}
Let $R$ be a commutative ring such that $\bbQ \subseteq R$.
Suppose that $\calh_?^*$ is a proper
cohomology theory with values in $R$-modules which comes with a
multiplicative structure. Suppose that the $R\Sub(G;\calf)$-module
$\calh^q_G(G/?)$ of \eqref{Sub(G,F)-coefficients H_G(G/?))}, 
which sends $G/H$ to $\calh^q_G(G/H)$, is injective for each $q \in
\bbZ$.

Then the natural transformation of equivariant cohomology theories appearing in
Theorem~\ref{the: equivariant Chern character and injective coeff}
$$
\ch^*_? \colon \calh^*_?   ~ \to ~ \calbh^*_?
$$ 
is compatible with the given multiplicative structure on $\calh^*_?$
and the induced multiplicative structure on $\calbh^*_?$.
\end{theorem}

\begin{remark}[External products and restriction structures] 
\label{rem: External products and restriction structures} \em
One can also define an \emph{external product} for a proper equivariant cohomology theory
$\calh^*_?$ with values in $R$-modules. It assigns to every two groups $G$ and $H$,
a proper $G$-$CW$-pair $(X,A)$, a proper $H$-$CW$-pair´$(Y,B)$ and $p,q \in \bbZ$ an
$R$-homomorphism
$$\times \colon \calh^p_G(X,A) \otimes_R \calh^q_H(Y,B) ~ \to ~ 
\calh^{p+q}_{G \times H}((X,A) \times (Y,B)).$$
One requires graded commutativity, associativity, the existence of a unit
$1 \in \calh^0_{\{1\}}(\pt)$ and compatibility with the induction structure and the
boundary homomorphism associated to a pair. One can show that $\calbh^*_?$ inherits
an external product and prove the analogon of 
Theorem~\ref{the: equivariant Chern character and multiplicative structures}
for external products.

One can also introduce the notion of a \emph{restriction structure} on $\calh^*_?$.
It yields for every injective group homomorphism $\alpha \colon H \to G$, every
proper $G$-$CW$-pair $(X,A)$ and $p \in \bbZ$ an $R$-homomorphism
$$\res_{\alpha} \colon \calh^p_G(X,A)) ~ \to ~ \calh^p_H(\res_{\alpha}(X,A)).$$
Again certain axioms are required such as compatibility with the boundary homomorphism
associated to pair, compatibility with induction for group isomorphisms
$\alpha \colon H \xrightarrow{\cong} G$, compatibility with conjugation, 
the double coset formula and compatibility for projections onto quotients under free actions.
One can show that $\calbh^*_?$ inherits
a restriction structure and prove the analogon of 
Theorem~\ref{the: equivariant Chern character and multiplicative structures}
for restriction structures.

An external product together with a restriction structure yields a multiplicative
structure as follows. Consider 
$G$-$CW$-pairs $(X,A)$ and $(X,B)$. Let $d \colon G \to G \times G$ and
$D \colon (X,A \cup B) \to (X,A) \times (X,B)$ be the diagonal maps.
Define 
\begin{eqnarray}
\cup \colon \calh^m_G(X,A) \otimes_R \calh^n_G(X,A) & \to &
\calh^{m+n}_G(X,A \cup B)
\label{cup product} 
\end{eqnarray}
to be the composite
\begin{multline*}
\calh^m_G(X,A) \otimes_R \calh^n_G(X,A) 
~ \xrightarrow{\times} \calh^{m+n}_{G \times G}((X,A) \times (X,B))
\\
~ \xrightarrow{\res_d} \calh^{m+n}_{G}((X,A) \times (X,B))
~ \xrightarrow{\calh^{m+n}_{G}(D)} \calh^{m+n}_{G}(X,A \cup B).
\end{multline*}
\em
\end{remark}


\bibliographystyle{abbrv}
\bibliography{dbdef,dbpub,dbpre,dbecextra}

\end{document}